\documentclass{amsart}

\usepackage{lipsum}
\usepackage{amsfonts}
\usepackage{graphicx}
\usepackage{epstopdf}
\usepackage{algorithm}
\ifpdf%
  \DeclareGraphicsExtensions{.eps,.pdf,.png,.jpg}
\else
  \DeclareGraphicsExtensions{.eps}
\fi
\usepackage{amsopn}
\usepackage{booktabs}

\usepackage[utf8]{inputenc}
\usepackage{etoolbox}
\usepackage{pgfplots}
\usepackage[english]{babel}
\usepackage{amsmath}
\usepackage{amssymb}
\usepackage{amsthm}
\usepackage{mathtools}
\DeclareMathAlphabet{\mathsfit}{T1}{\sfdefault}{\mddefault}{\sldefault}
\usepackage[sort, numbers]{natbib}
\usepackage{tikz}
\usepackage{diagbox}
\usepackage{widebar}
\usepackage{uniinput}
\usepackage{hyperref}
\usepackage[capitalize]{cleveref}
\usepackage[noend]{algpseudocode}
\usepackage{subfig}
\usepackage{xfrac}
\usepackage[flushleft]{threeparttable}
\usepackage{setspace}
\usepackage{enumitem}
\usepackage{accents}
\usepackage{colortbl}
\usepackage{scalerel}
\usepackage{multirow}

\usepackage{xspace}
\newcommand*{\ie}{i.e.\@\xspace}%
\newcommand*{\eg}{e.g.\@\xspace}%
\robustify{\subref}
\robustify{\widebar}
\robustify{\phantom}

\newlist{loops}{enumerate}{1}
\setlist[loops,1]
{
  label=(L\arabic*),
  ref=L\arabic*
}

\newlist{assignments}{enumerate}{1}
\setlist[assignments,1]
{
  label=(A\arabic*),
  ref=A\arabic*
}

\crefname{equation}{}{}
\crefname{loopsi}{}{}
\creflabelformat{loopsi}{(#2#1#3)}
\crefname{assignmentsi}{}{}
\creflabelformat{assignmentsi}{(#2#1#3)}

\makeatletter
\newcommand*{\algrule}[1][\algorithmicindent]{%
  \makebox[#1][l]{%
    \hspace*{.2em}
    \vrule height .75\baselineskip depth .25\baselineskip
  }
}

\newcount\ALG@printindent@tempcnta
\def\ALG@printindent{%
    \ifnum \theALG@nested>0
    \ifx\ALG@text\ALG@x@notext
    \else
    \unskip
    \ALG@printindent@tempcnta=1
    \loop
    \algrule[\csname ALG@ind@\the\ALG@printindent@tempcnta\endcsname]%
    \advance \ALG@printindent@tempcnta 1
    \ifnum \ALG@printindent@tempcnta<\numexpr\theALG@nested+1\relax
    \repeat
    \fi
    \fi
}
\patchcmd{\ALG@doentity}{\noindent\hskip\ALG@tlm}{\ALG@printindent}{}{\errmessage{failed to patch}}
\patchcmd{\ALG@doentity}{\item[]\nointerlineskip}{}{}{} 
\makeatother

\usepackage{tikz-3dplot}
\usetikzlibrary{shapes.geometric}
\usetikzlibrary{intersections}
\usetikzlibrary{matrix,arrows,positioning,calc}
\usetikzlibrary{fit}
\usepackage{IMjournal}

\newtheorem{theorem}{Theorem}[section]
\newtheorem{corollary}[theorem]{Corollary}
\theoremstyle{remark}
\newtheorem{remark}[theorem]{Remark}
\newtheorem{example}[theorem]{Example}
\numberwithin{equation}{section}

\newcommand{\slowGraphics}[1]{%
  \ifx\OMITSLOWGRAPHICS\undefined{#1}%
  \else{\centering\footnotesize Graphics omited!}\fi%
}%

\newcommand{\myTitle}{$k\ell$-refinement: An adaptive mesh refinement scheme for hiearchical hybrid grids}

\newcommand\hyteg{\textsc{{HyTeG}}\@\xspace}

\renewcommand\tilde[1]{\widetilde{#1}}

\newcommand\uh{\mathbf{u}}

\newcommand\wh{\mathbf{w}}
\newcommand\bh{\mathbf{b}}

\newcommand\eest{\tilde{\mathbf{e}}}

\newcommand \xx{{\bf x}}

\renewcommand\L{\mathsfit{L}}
\newcommand\K{\mathsfit{K}}
\newcommand\hmean{\bar{h}}
\newcommand\Hmean{\bar{H}}

\newcommand\Tri{\mathcal{T}}

\newcommand\Triloc{\Tri_\mathrm{loc}}

\newcommand\SET[1]{\left\{#1\right\}}
\newcommand\SEQ[2]{\left(#1\right)_{#2}}
\newcommand\setsep{\,\big|\,}

\newcommand\Pspace[1]{\mathbb{P}_#1}

\newcommand\Pj{\Pspace{1}}

\newcommand\lowBar[2]{\raisebox{-0.2ex}{$#1|$}_{#2}}
\newcommand\restrict[1]{\mathpalette\lowBar{#1}}
\newcommand\norm[2]{\left\|#1\right\|_{#2}}

\newcommand\transposesign{\hspace{-0.6pt}\mathsf{T}}
\renewcommand\t{^{\transposesign}}

\newcommand\Oh[1]{\bigO(\hmean^{#1})}

\newcommand\errest[1]{\tilde{e}_{#1}}
\newcommand\conv{θ}
\newcommand\convest[1]{\tilde{\conv}_{#1}}
\newcommand\convmin{\conv_{*}}
\newcommand\convmax{\conv_{\mathsf{x}}}
\newcommand\errratio{ϱ_{\Hmean}^{\hmean}}

\newcommand\Rc{\mathcal{R}}
\newcommand\Gc{\mathcal{G}}

\newcommand\bigO{\mathcal{O}}

\DeclarePairedDelimiter\ceil{\lceil}{\rceil}

\usepackage[textsize=footnotesize]{todonotes}
\ifpdf
\hypersetup{
  pdftitle={\myTitle},
  pdfauthor={B. Mann}
}
\fi

\graphicspath{{./figures/}}

\begin{document}

  \title{$k\ell$-refinement: An adaptive mesh refinement scheme for hierarchical hybrid grids}
  \author{%
    \|Benjamin |Mann|, %
    Chair for System Simulation, Friedrich-Alexander-University,  Erlangen-Nürnberg, Germany,\\%
    \|Ulrich |Rüde|, %
    Department of Applied Mathematics, VSB-Technical University of Ostrava, Ostrava, Czech Republic and CERFACS, Toulouse, France%
  }%



  \abstract
    This work introduces an adaptive mesh refinement technique for hierarchical hybrid grids
with the goal to reach scalability and maintain excellent performance on massively parallel computer systems.
On the block structured hierarchical hybrid grids, this is accomplished by
using classical, unstructured refinement only on the coarsest level of the hierarchy, while keeping the number of
structured refinement levels constant on the whole domain.
This leads to a compromise where the excellent performance characteristics
of hierarchical hybrid grids can be maintained at the price that the flexibility of generating locally refined meshes is constrained.
Furthermore, mesh adaptivity
often relies on a posteriori error estimators or error indicators
that tend to become computationally expensive.
Again with the goal of preserving scalability and performance,
a method is proposed that leverages the grid hierarchy and the full multigrid scheme that
generates a natural sequence of approximations on the nested hierarchy of grids.
This permits to compute a cheap error estimator that is well-suited for large-scale parallel computing.
We present the theoretical foundations for both global and local error estimates and present a rigorous analysis of their effectivity.
The proposed method, including error estimator and the adaptive coarse grid refinement,
is implemented in the finite element framework \hyteg.
Extensive numerical experiments are conducted to validate the effectiveness, as well as performance and scalability.
  \endabstract



  \section{Introduction}\label{sec:intro}
Modern parallel supercomputers enable numerical simulations 
of enormous size so that systems with extremely large linear systems must be solved.
One guiding example for the present work is Earth's mantle convection, where the number of degrees of freedom (DoF) can reach up to $10^{13}$~\cite{GMEINER2016509}.
Even a single vector of this size is approximately 80 TByte large.

Efficiently solving these massive systems requires appropriate algorithms and scalable implementations.
The most crucial requirement for the linear solver is optimal order complexity, achievable through multigrid (MG) methods or MG-preconditioned Krylov solvers, see \eg~\cite{elman2001multigrid}.
Additionally, using matrix-free methods~\cite{KRONBICHLER} may be necessary.
Multigrid methods are most effective on meshes with an inherent hierarchical structure, such as block-structured grids resulting from successive regular refinement of an initial, possibly unstructured, coarse grid.
The free open source framework \hyteg\footnote{\url{https://i10git.cs.fau.de/hyteg/hyteg}} is based on such a hybrid mesh construction~\cite{hyteg, kohl2024fundamental}.
This approach does not only facilitate the straightforward design of a geometric MG scheme but also inherits additional advantages of structured grids.
Notably, this includes computational kernels typical for structured grids that can inherently be well optimized for modern processor architectures. This can lead to significantly faster execution times~\cite{mayr2022noninvasive, boehm2024}.

Many real-world applications require adapting the grid resolution non-uniformly in order to compute
a sufficiently accurate solution with a given number of grid points.
Thus, adaptive mesh refinement (AMR) is an important feature of many modern simulation frameworks.
Achieving mesh adaptivity is in principle straightforward when using fully unstructured grids.
On block-structured grids, however, there are certain restrictions regarding the choice of refinement strategy.
A widely used approach is the following:
Instead of using the same level of structured refinement everywhere in the domain,
the refinement depth is chosen individually on each block, \eg~\cite{brandt1977multi,mccormick1986, mccorquodale2004node, gradl, zhang2019amrex, munch2022efficient}.
One downside of this strategy is that it results in hanging nodes that must be treated appropriately.
Preserving the optimal order complexity of the resulting multigrid algorithms also requires special care~\cite{munch2022efficient}.
Most importantly though,  the performance advantage of block-structured grids cannot be fully exploited
when the maximal level of structured refinement is only used in a few blocks covering the computational domain.

In this work, we present an alternative approach, both maintaining the conformity of the resulting mesh and keeping the number of local refinement levels constant across the domain.
This is achieved by applying a specially modified, adaptive refinement scheme exclusively to the coarse grid.
This strategy allows leveraging the efficiency of the algorithms and data structures implemented in \hyteg, as well as the automatically generated, highly optimized compute kernels.

In addition to refinement strategies, appropriate refinement criteria are essential.
Typically, local a posteriori error indicators are used to determine where the mesh should be refined.
It is also often necessary to evaluate the global accuracy, especially when AMR is performed iteratively until reaching sufficiently small discretization errors.
There are several well-established methods for estimating local and/or global errors.
Prominent classical examples include gradient recovery estimators such as the Zienkiewicz and Zhu estimator~\cite{zizu} and element residual methods~\cite{babuska1978}.
For a comprehensive overview, see~\cite{verfurth1996,babuska2001reliability,ainsworth1997}.

Many of the more recent 
error estimators are computationally expensive as they involve solving additional subproblems, see \eg~\cite{luce2004local,richter2015goal,papez2020sharp}.
In this work, we propose an error estimator, tailored to our data structures and algorithms that makes 
use of the data that is inherently computed by the full multigrid (FMG) solver \cite{brandt1977multi, trottenberg2000multigrid}.
In particular, approximate solutions from coarser grid levels are utilized, which are computed without extra cost within the FMG scheme, see e.g.\ \cite{bai1987local}.
The fine-grid solution then serves as a proxy for the 
correct solution.
This concept is akin to traditional predictor-corrector strategies for ordinary differential equations~\cite{hairer1993solving}.
Similar methods for partial differential equations (PDEs) have been explored in \eg~\cite{AULISA2018224,Ferraz-Leite2010},
though primarily in the context of unstructured grids.
The hierarchical grid structure can be leveraged both in the theoretical analysis and in the practical implementation within an HPC framework.
The resulting error estimator is problem-agnostic, computationally inexpensive, highly scalable, and straightforward to implement.

  \section{Hybrid tetrahedral grids}\label{sec:hhg}
\begin{figure}
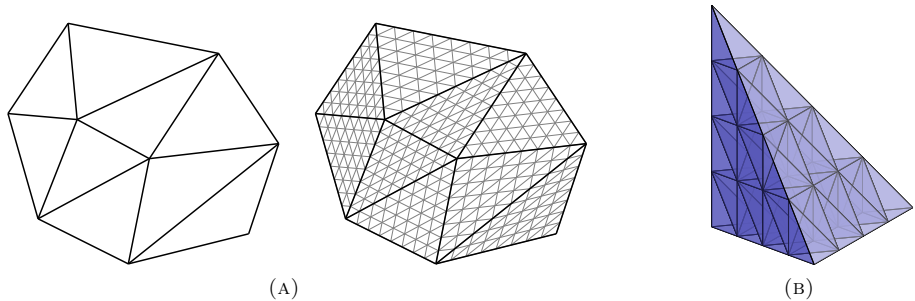

    \centering
    \subfloat[\label{fig:hhg:2d}]
    {
        \slowGraphics{
            \tdplotsetmaincoords{0}{0}
\begin{tikzpicture}[tdplot_main_coords, scale=0.4]

\coordinate (N0) at (0.0, 0.0);
\coordinate (N1) at (2.0, 3.0);
\coordinate (N2) at (4.0, -5.0);
\coordinate (N3) at (7.0, 2.0);
\coordinate (N4) at (8.0, -4.0);
\coordinate (N5) at (9.0, -1.0);
\coordinate (N6) at (2.3, -0.19999999999999996);
\coordinate (N7) at (1.0, -3.5);
\coordinate (N8) at (4.7, -1.5);
\draw [line width=0.2mm, color=black] (N0) -- (N1);
\draw [line width=0.2mm, color=black] (N0) -- (N6);
\draw [line width=0.2mm, color=black] (N0) -- (N7);
\draw [line width=0.2mm, color=black] (N1) -- (N3);
\draw [line width=0.2mm, color=black] (N1) -- (N6);
\draw [line width=0.2mm, color=black] (N2) -- (N4);
\draw [line width=0.2mm, color=black] (N2) -- (N5);
\draw [line width=0.2mm, color=black] (N2) -- (N7);
\draw [line width=0.2mm, color=black] (N2) -- (N8);
\draw [line width=0.2mm, color=black] (N3) -- (N5);
\draw [line width=0.2mm, color=black] (N3) -- (N6);
\draw [line width=0.2mm, color=black] (N3) -- (N8);
\draw [line width=0.2mm, color=black] (N4) -- (N5);
\draw [line width=0.2mm, color=black] (N5) -- (N8);
\draw [line width=0.2mm, color=black] (N6) -- (N7);
\draw [line width=0.2mm, color=black] (N6) -- (N8);
\draw [line width=0.2mm, color=black] (N7) -- (N8);

\end{tikzpicture}
            \input{figures/hhg2d_lvl3}
        }
    }
    \hfill
    \subfloat[\label{fig:hhg:3d}]
    {
        \slowGraphics{
            \input{figures/tet_lvl2}
        }
    }
    \caption{
        \subref{fig:hhg:2d}:~Macro grid $\Tri_0$ of a two dimensional domain and corresponding fine grid $\Tri_{3}$.
        \subref{fig:hhg:3d}:~Single macro tetrahedron on local refinement level $\ell=2$.
    }
    \label{fig:hhg}
\end{figure}

In the scope of this work, we consider block-structured, simplicial meshes often referred to as \emph{Hierarchical Hybrid Grids} (HHG)~\cite{hhg}.
These meshes are constructed based on an unstructured but rather coarse triangulation $\Tri_0$, which is then uniformly refined $\L$ times following the refinement algorithm of Bey~\cite{bey}, see also \cite{blaheta2003nested, axelsson2004two, blaheta2006hierarchical},
where a strengthened Cauchy-Schwarz inequality for such a construction is proven as the basis for fast multigrid convergence.
A sketch illustrating the refinement procedure is displayed in \cref{fig:hhg}.
The mesh $\Tri_0$ and its elements $T∈\Tri_0$ are referred to as \emph{macro grid} and \emph{macro elements}, respectively.
These HHG have several advantages compared to fully unstructured grids, while still preserving some of their flexibility.
Most importantly, the structuredness allows for various performance optimizations as well as the automated generation of compute kernels.
Clearly, the impact of said optimizations scales with the level of regular refinement.
Therefore, $\L$ should be chosen as large as possible, or, equivalently, $\Tri_0$ should be as coarse as possible.
In our in-house code \hyteg, the coarse grid is used as basis for parallel data structures.
Each process then only stores the data corresponding to a subset $\Triloc\subset\Tri_0$, which means that $\L$ is only bounded by the available memory per process.
Apart from performance aspects, the HHG structure can also be exploited algorithmically.
Most notably, the sequence of nested triangulations $\SEQ{\Tri_\ell}{\ell=0,\dots,\L}$,
leads to a straightforward 
construction of geometric multigrid methods.
In combination, these properties allow for very efficient and highly scalable implementations of matrix-free operators and algorithms in \hyteg~\cite{hyteg, kohl2022, boehm2024}.

  \section{$k\ell$-refinement}\label{sec:rg}

An apparent downside of HHG is its limited flexibility with respect to adaptive refinement.
While it seems straightforward to use different levels of refinement on different macro elements,
this has several disadvantages.
Such a construction leads to hanging nodes and potential load imbalance in parallel computing.
Additionally, a loss of performance must be expected when the level of structured meshes is reduced,
the performance benefit of using HHG might be compromised.

In this work, we therefore employ an alternative strategy that preserves both the conformity of the 
meshes and the uniformity of the local refinement.
The latter is achieved by restricting the adaptive 
refinement process to the coarse grid $\Tri_0$, while keeping the structured local refinement as described in \cref{sec:hhg}.
The mesh corresponding to HHG level $\ell$ of a $k$ times adaptively refined coarse grid will be denoted $\Tri_\ell^k$, omitting either index when there is no ambiguity.
Canonically, AMR on the coarse grid, and uniform refinement are referred to as $k$- and $\ell$-refinement, respectively.
Note that there are two natural ways to obtain the grid $\Tri_\ell^k$:
One can either apply $k$-refinement to $\Tri^0$ and then use $\Tri^k$ as basis for HHG,
or apply $k$-refinement to the coarse grid $\Tri_0$ of an HHG hierarchy $\SEQ{\Tri_\ell}{\ell}$.
We refer to the latter as $k\ell$-refinement, as both strategies are intertwined here.
The two approaches are summarized in \cref{alg:kplusl,alg:kl}, respectively.

  \begin{algorithm}
    \footnotesize
    \begin{algorithmic}[1]
      \setstretch{1.5}
      \State Given a coarse grid $\Tri_0^0$
      \For{$k = 0,...,\K-1$}
        \State Solve the problem on $\Tri_0^{k}$.
        \State Estimate error on each macro element.
        \State Obtain $\Tri_0^{k+1}$ by AMR of $\Tri_0^{k}$.
      \EndFor
      \For{$\ell=0,...,\L-1$}
        \State Obtain $\Tri_{\ell+1}^{\K}$ by uniform refinement of $\Tri_{\ell-1}^{\K}$.
      \EndFor
      \State Solve the problem on $\Tri_\L^{\K}$.
    \end{algorithmic}
    \caption{$k$-Refinement + $\ell$-Refinement}
    \label{alg:kplusl}
  \end{algorithm}
  \begin{algorithm}
    \footnotesize
    \begin{algorithmic}[1]
      \setstretch{1.5}
      \State Given a coarse grid $\Tri_0^0$
      \For{$k = 0,...,\K-1$}
        \For{$\ell=0,...,\L-1$}
          \State Obtain $\Tri_{\ell+1}^{k}$ by uniform refinement of $\Tri_{\ell-1}^{k}$.
        \EndFor
        \State Solve the problem on $\Tri_{\L}^{k}$.
        \State Estimate error on each macro element.
        \State Obtain $\Tri_0^{k+1}$ by AMR of $\Tri_0^{k}$.
      \EndFor
      \State Solve the problem on $\Tri_\L^{\K}$.
    \end{algorithmic}
    \caption{$k\ell$-Refinement}
    \label{alg:kl}
  \end{algorithm}

To preserve the conformity of the mesh, an appropriate algorithm such as \eg newest node
bisection~\cite{stals, stals1995adaptive, fung2024fault} or red-green refinement (RG)~\cite{bey, pltmg},
must be employed for $k$-refinement.
In the scope of this work, RG is considered.
In the following, the key aspects of the algorithm are highlighted.
A detailed description can be found in \cite{bey}.
%

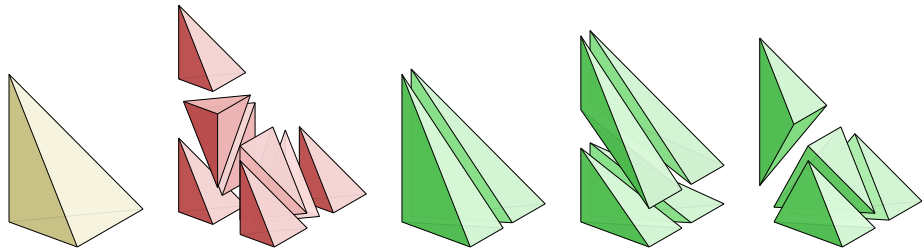
\begin{figure}
  \centering
  \slowGraphics{
    \tdplotsetmaincoords{0}{0}
\begin{tikzpicture}[tdplot_main_coords, scale=2]

\coordinate (A0) at (0.8910065241883679, 0.25202731045750026, 0.4459490798461028);
\coordinate (B0) at (0.0, 0.1669579738068738, 0.0);
\coordinate (C0) at (0.0, 1.1492452245355624, -0.18738131458572482);
\coordinate (D0) at (0.45399049973954675, 0.0, -0.8752243490263167);
\fill[black!46!yellow!84, opacity=0.8] (A0) -- (B0) -- (C0) -- cycle;
\fill[black!46!yellow!84, opacity=0.8] (A0) -- (B0) -- (D0) -- cycle;
\draw [line width=0.05mm, color=black, solid] (A0) -- (B0);
\fill[black!37!yellow!68, opacity=0.9] (B0) -- (C0) -- (D0) -- cycle;
\fill[black!6!yellow!11, opacity=0.9] (A0) -- (C0) -- (D0) -- cycle;
\draw [line width=0.05mm, color=black] (B0) -- (C0);
\draw [line width=0.05mm, color=black] (A0) -- (D0);
\draw [line width=0.05mm, color=black] (C0) -- (D0);
\draw [line width=0.05mm, color=black] (A0) -- (C0);
\draw [line width=0.05mm, color=black] (B0) -- (D0);

\end{tikzpicture}
    \tdplotsetmaincoords{0}{0}
\begin{tikzpicture}[tdplot_main_coords, scale=2]

\coordinate (A0) at (0.44550326209418395, 0.2762758316549365, 0.1798085790594357);
\coordinate (B0) at (0.0, 0.2337411633296233, -0.04316596086361572);
\coordinate (C0) at (0.0, 0.7248847886939676, -0.13685661815647812);
\coordinate (D0) at (0.22699524986977337, 0.15026217642618642, -0.4807781353767741);
\fill[black!46!red!84, opacity=0.8] (A0) -- (B0) -- (C0) -- cycle;
\fill[black!46!red!84, opacity=0.8] (A0) -- (B0) -- (D0) -- cycle;
\draw [line width=0.05mm, color=black, solid] (A0) -- (B0);
\fill[black!37!red!68, opacity=0.9] (B0) -- (C0) -- (D0) -- cycle;
\fill[black!6!red!11, opacity=0.9] (A0) -- (C0) -- (D0) -- cycle;
\draw [line width=0.05mm, color=black] (B0) -- (C0);
\draw [line width=0.05mm, color=black] (A0) -- (D0);
\draw [line width=0.05mm, color=black] (C0) -- (D0);
\draw [line width=0.05mm, color=black] (A0) -- (C0);
\draw [line width=0.05mm, color=black] (B0) -- (D0);

\coordinate (A1) at (0.5078737187873696, 0.4688652628589312, 0.17542256487737518);
\coordinate (B1) at (0.5078737187873696, 0.9600088882232756, 0.08173190758451276);
\coordinate (C1) at (0.06237045669318575, 0.9174742198979624, -0.14124263233853865);
\coordinate (D1) at (0.2893657065629591, 0.3428516076301811, -0.48516414955883463);
\fill[black!46!red!84, opacity=0.8] (A1) -- (B1) -- (C1) -- cycle;
\fill[black!53!red!97, opacity=0.8] (A1) -- (C1) -- (D1) -- cycle;
\draw [line width=0.05mm, color=black, solid] (A1) -- (C1);
\fill[black!13!red!23, opacity=0.9] (A1) -- (B1) -- (D1) -- cycle;
\fill[black!6!red!11, opacity=0.9] (B1) -- (C1) -- (D1) -- cycle;
\draw [line width=0.05mm, color=black] (A1) -- (B1);
\draw [line width=0.05mm, color=black] (B1) -- (C1);
\draw [line width=0.05mm, color=black] (A1) -- (D1);
\draw [line width=0.05mm, color=black] (C1) -- (D1);
\draw [line width=0.05mm, color=black] (B1) -- (D1);

\coordinate (A2) at (0.4772825970759523, 1.0111270840422588, -0.023866924457557337);
\coordinate (B2) at (0.03177933498176827, 0.9685924157169454, -0.24684146438060875);
\coordinate (C2) at (0.25877458485154164, 0.39396980344916427, -0.5907629816009047);
\coordinate (D2) at (0.25877458485154164, 0.8851134288135085, -0.6844536388937671);
\fill[black!53!red!97, opacity=0.8] (A2) -- (B2) -- (C2) -- cycle;
\fill[black!13!red!24, opacity=0.9] (A2) -- (B2) -- (D2) -- cycle;
\fill[black!37!red!68, opacity=0.9] (B2) -- (C2) -- (D2) -- cycle;
\fill[black!13!red!23, opacity=0.9] (A2) -- (C2) -- (D2) -- cycle;
\draw [line width=0.05mm, color=black] (A2) -- (B2);
\draw [line width=0.05mm, color=black] (B2) -- (C2);
\draw [line width=0.05mm, color=black] (A2) -- (D2);
\draw [line width=0.05mm, color=black] (C2) -- (D2);
\draw [line width=0.05mm, color=black] (A2) -- (C2);
\draw [line width=0.05mm, color=black] (B2) -- (D2);

\coordinate (A3) at (0.44550326209418395, 1.1603343573107563, 0.011165395932283356);
\coordinate (B3) at (0.0, 1.1177996889854431, -0.21180914399076806);
\coordinate (C3) at (0.0, 1.6089433143497873, -0.30549980128363047);
\coordinate (D3) at (0.22699524986977337, 1.0343207020820062, -0.6494213185039264);
\fill[black!46!red!84, opacity=0.8] (A3) -- (B3) -- (C3) -- cycle;
\fill[black!46!red!84, opacity=0.8] (A3) -- (B3) -- (D3) -- cycle;
\draw [line width=0.05mm, color=black, solid] (A3) -- (B3);
\fill[black!37!red!68, opacity=0.9] (B3) -- (C3) -- (D3) -- cycle;
\fill[black!6!red!11, opacity=0.9] (A3) -- (C3) -- (D3) -- cycle;
\draw [line width=0.05mm, color=black] (B3) -- (C3);
\draw [line width=0.05mm, color=black] (A3) -- (D3);
\draw [line width=0.05mm, color=black] (C3) -- (D3);
\draw [line width=0.05mm, color=black] (A3) -- (C3);
\draw [line width=0.05mm, color=black] (B3) -- (D3);

\coordinate (A4) at (1.247409133863715, 0.35283823464050035, 0.5811627509209283);
\coordinate (B4) at (0.8019058717695311, 0.31030356631518713, 0.3581882109978768);
\coordinate (C4) at (0.8019058717695311, 0.8014471916795316, 0.2644975537050144);
\coordinate (D4) at (1.0289011216393045, 0.22682457941175022, -0.07942396351528155);
\fill[black!46!red!84, opacity=0.8] (A4) -- (B4) -- (C4) -- cycle;
\fill[black!46!red!84, opacity=0.8] (A4) -- (B4) -- (D4) -- cycle;
\draw [line width=0.05mm, color=black, solid] (A4) -- (B4);
\fill[black!37!red!68, opacity=0.9] (B4) -- (C4) -- (D4) -- cycle;
\fill[black!6!red!11, opacity=0.9] (A4) -- (C4) -- (D4) -- cycle;
\draw [line width=0.05mm, color=black] (B4) -- (C4);
\draw [line width=0.05mm, color=black] (A4) -- (D4);
\draw [line width=0.05mm, color=black] (C4) -- (D4);
\draw [line width=0.05mm, color=black] (A4) -- (C4);
\draw [line width=0.05mm, color=black] (B4) -- (D4);

\coordinate (A5) at (0.7089442933649279, 0.28670680101761825, 0.23448963538758028);
\coordinate (B5) at (0.7089442933649279, 0.7778504263819626, 0.14079897809471786);
\coordinate (C5) at (0.9359395432347011, 0.20322781411418134, -0.2031225391255781);
\coordinate (D5) at (0.49043628114051724, 0.16069314578886812, -0.42609707904862953);
\fill[black!22!red!40, opacity=0.8] (A5) -- (B5) -- (C5) -- cycle;
\fill[black!47!red!85, opacity=0.8] (A5) -- (B5) -- (D5) -- cycle;
\fill[black!46!red!84, opacity=0.8] (A5) -- (C5) -- (D5) -- cycle;
\draw [line width=0.05mm, color=black, solid] (A5) -- (B5);
\draw [line width=0.05mm, color=black, solid] (A5) -- (C5);
\draw [line width=0.05mm, color=black, solid] (A5) -- (D5);
\fill[black!6!red!11, opacity=0.9] (B5) -- (C5) -- (D5) -- cycle;
\draw [line width=0.05mm, color=black] (C5) -- (D5);
\draw [line width=0.05mm, color=black] (B5) -- (C5);
\draw [line width=0.05mm, color=black] (B5) -- (D5);

\coordinate (A6) at (0.6259112487196519, 0.7987253449460459, -0.12334066541204267);
\coordinate (B6) at (0.8529064985894252, 0.22410273267826453, -0.46726218263233865);
\coordinate (C6) at (0.4074032364952413, 0.1815680643529513, -0.69023672255539);
\coordinate (D6) at (0.4074032364952413, 0.6727116897172956, -0.7839273798482524);
\fill[black!53!red!97, opacity=0.8] (A6) -- (B6) -- (C6) -- cycle;
\fill[black!47!red!85, opacity=0.8] (A6) -- (C6) -- (D6) -- cycle;
\draw [line width=0.05mm, color=black, solid] (A6) -- (C6);
\fill[black!6!red!11, opacity=0.9] (A6) -- (B6) -- (D6) -- cycle;
\fill[black!13!red!24, opacity=0.9] (B6) -- (C6) -- (D6) -- cycle;
\draw [line width=0.05mm, color=black] (A6) -- (B6);
\draw [line width=0.05mm, color=black] (B6) -- (C6);
\draw [line width=0.05mm, color=black] (A6) -- (D6);
\draw [line width=0.05mm, color=black] (C6) -- (D6);
\draw [line width=0.05mm, color=black] (B6) -- (D6);

\coordinate (A7) at (0.8540947118597759, 0.1260136552287501, -0.6078933350642495);
\coordinate (B7) at (0.4085914497655921, 0.08347898690343689, -0.8308678749873009);
\coordinate (C7) at (0.4085914497655921, 0.5746226122677811, -0.9245585322801633);
\coordinate (D7) at (0.6355866996353654, 0.0, -1.2684800495004591);
\fill[black!46!red!84, opacity=0.8] (A7) -- (B7) -- (C7) -- cycle;
\fill[black!46!red!84, opacity=0.8] (A7) -- (B7) -- (D7) -- cycle;
\draw [line width=0.05mm, color=black, solid] (A7) -- (B7);
\fill[black!37!red!68, opacity=0.9] (B7) -- (C7) -- (D7) -- cycle;
\fill[black!6!red!11, opacity=0.9] (A7) -- (C7) -- (D7) -- cycle;
\draw [line width=0.05mm, color=black] (B7) -- (C7);
\draw [line width=0.05mm, color=black] (A7) -- (D7);
\draw [line width=0.05mm, color=black] (C7) -- (D7);
\draw [line width=0.05mm, color=black] (A7) -- (C7);
\draw [line width=0.05mm, color=black] (B7) -- (D7);

\end{tikzpicture}
    \tdplotsetmaincoords{0}{0}
\begin{tikzpicture}[tdplot_main_coords, scale=2]

\coordinate (A0) at (0.9521887676112029, 0.2873111339215503, 0.45209929813994076);
\coordinate (B0) at (0.061182243422834987, 0.20224179727092384, 0.006150218293837946);
\coordinate (C0) at (0.061182243422834987, 1.1845290479996124, -0.18123109629188688);
\coordinate (D0) at (0.7336807553867923, 0.16129747869280014, -0.20848741629626902);
\fill[black!46!green!84, opacity=0.8] (A0) -- (B0) -- (C0) -- cycle;
\fill[black!46!green!84, opacity=0.8] (A0) -- (B0) -- (D0) -- cycle;
\draw [line width=0.05mm, color=black, solid] (A0) -- (B0);
\fill[black!24!green!43, opacity=0.9] (B0) -- (C0) -- (D0) -- cycle;
\fill[black!6!green!11, opacity=0.9] (A0) -- (C0) -- (D0) -- cycle;
\draw [line width=0.05mm, color=black] (B0) -- (C0);
\draw [line width=0.05mm, color=black] (A0) -- (D0);
\draw [line width=0.05mm, color=black] (C0) -- (D0);
\draw [line width=0.05mm, color=black] (A0) -- (C0);
\draw [line width=0.05mm, color=black] (B0) -- (D0);

\coordinate (A1) at (0.0, 0.1669579738068738, -0.1788140617483008);
\coordinate (B1) at (0.0, 1.1492452245355624, -0.36619537633402566);
\coordinate (C1) at (0.6724985119639573, 0.12601365522875013, -0.39345169633840776);
\coordinate (D1) at (0.4539904997395467, 0.0, -1.0540384107746175);
\fill[black!36!green!65, opacity=0.8] (A1) -- (B1) -- (C1) -- cycle;
\fill[black!46!green!84, opacity=0.8] (A1) -- (C1) -- (D1) -- cycle;
\draw [line width=0.05mm, color=black, solid] (A1) -- (C1);
\fill[black!37!green!68, opacity=0.9] (A1) -- (B1) -- (D1) -- cycle;
\fill[black!6!green!11, opacity=0.9] (B1) -- (C1) -- (D1) -- cycle;
\draw [line width=0.05mm, color=black] (A1) -- (B1);
\draw [line width=0.05mm, color=black] (B1) -- (C1);
\draw [line width=0.05mm, color=black] (A1) -- (D1);
\draw [line width=0.05mm, color=black] (C1) -- (D1);
\draw [line width=0.05mm, color=black] (B1) -- (D1);

\end{tikzpicture}
    \tdplotsetmaincoords{0}{0}
\begin{tikzpicture}[tdplot_main_coords, scale=2]

\coordinate (A0) at (0.9521887676112029, 0.2873111339215503, 0.4652159901609415);
\coordinate (B0) at (0.061182243422834987, 0.20224179727092384, 0.019266910314838685);
\coordinate (C0) at (0.061182243422834987, 0.6933854226352681, -0.07442374697802373);
\coordinate (D0) at (0.7336807553867923, 0.16129747869280014, -0.19537072427526828);
\fill[black!46!green!84, opacity=0.8] (A0) -- (B0) -- (C0) -- cycle;
\fill[black!46!green!84, opacity=0.8] (A0) -- (B0) -- (D0) -- cycle;
\draw [line width=0.05mm, color=black, solid] (A0) -- (B0);
\fill[black!24!green!43, opacity=0.9] (B0) -- (C0) -- (D0) -- cycle;
\fill[black!6!green!11, opacity=0.9] (A0) -- (C0) -- (D0) -- cycle;
\draw [line width=0.05mm, color=black] (B0) -- (C0);
\draw [line width=0.05mm, color=black] (A0) -- (D0);
\draw [line width=0.05mm, color=black] (C0) -- (D0);
\draw [line width=0.05mm, color=black] (A0) -- (C0);
\draw [line width=0.05mm, color=black] (B0) -- (D0);

\coordinate (A1) at (0.0, 0.1669579738068738, -0.16569736972730006);
\coordinate (B1) at (0.0, 0.6581015991712181, -0.2593880270201625);
\coordinate (C1) at (0.6724985119639573, 0.12601365522875013, -0.380335004317407);
\coordinate (D1) at (0.4539904997395467, 0.0, -1.0409217187536168);
\fill[black!36!green!65, opacity=0.8] (A1) -- (B1) -- (C1) -- cycle;
\fill[black!46!green!84, opacity=0.8] (A1) -- (C1) -- (D1) -- cycle;
\draw [line width=0.05mm, color=black, solid] (A1) -- (C1);
\fill[black!37!green!68, opacity=0.9] (A1) -- (B1) -- (D1) -- cycle;
\fill[black!6!green!11, opacity=0.9] (B1) -- (C1) -- (D1) -- cycle;
\draw [line width=0.05mm, color=black] (A1) -- (B1);
\draw [line width=0.05mm, color=black] (B1) -- (C1);
\draw [line width=0.05mm, color=black] (A1) -- (D1);
\draw [line width=0.05mm, color=black] (C1) -- (D1);
\draw [line width=0.05mm, color=black] (B1) -- (D1);

\coordinate (A2) at (0.9521887676112029, 0.5427058191110093, 0.4164968483686531);
\coordinate (B2) at (0.061182243422834987, 0.9487801078247272, -0.12314288877031218);
\coordinate (C2) at (0.061182243422834987, 1.4399237331890715, -0.2168335460631746);
\coordinate (D2) at (0.7336807553867923, 0.41669216388225916, -0.24408986606755673);
\fill[black!46!green!84, opacity=0.8] (A2) -- (B2) -- (C2) -- cycle;
\fill[black!53!green!97, opacity=0.8] (A2) -- (B2) -- (D2) -- cycle;
\draw [line width=0.05mm, color=black, solid] (A2) -- (B2);
\fill[black!24!green!43, opacity=0.9] (B2) -- (C2) -- (D2) -- cycle;
\fill[black!6!green!11, opacity=0.9] (A2) -- (C2) -- (D2) -- cycle;
\draw [line width=0.05mm, color=black] (B2) -- (C2);
\draw [line width=0.05mm, color=black] (A2) -- (D2);
\draw [line width=0.05mm, color=black] (C2) -- (D2);
\draw [line width=0.05mm, color=black] (A2) -- (C2);
\draw [line width=0.05mm, color=black] (B2) -- (D2);

\coordinate (A3) at (0.0, 0.9134962843606771, -0.3081071688124509);
\coordinate (B3) at (0.0, 1.4046399097250215, -0.40179782610531334);
\coordinate (C3) at (0.6724985119639573, 0.38140834041820915, -0.4290541461096955);
\coordinate (D3) at (0.4539904997395467, 0.255394685189459, -1.0896408605459054);
\fill[black!36!green!65, opacity=0.8] (A3) -- (B3) -- (C3) -- cycle;
\fill[black!53!green!97, opacity=0.8] (A3) -- (C3) -- (D3) -- cycle;
\draw [line width=0.05mm, color=black, solid] (A3) -- (C3);
\fill[black!37!green!68, opacity=0.9] (A3) -- (B3) -- (D3) -- cycle;
\fill[black!6!green!11, opacity=0.9] (B3) -- (C3) -- (D3) -- cycle;
\draw [line width=0.05mm, color=black] (A3) -- (B3);
\draw [line width=0.05mm, color=black] (B3) -- (C3);
\draw [line width=0.05mm, color=black] (A3) -- (D3);
\draw [line width=0.05mm, color=black] (C3) -- (D3);
\draw [line width=0.05mm, color=black] (B3) -- (D3);

\end{tikzpicture}
    \tdplotsetmaincoords{0}{0}
\begin{tikzpicture}[tdplot_main_coords, scale=2]

\coordinate (A0) at (0.44550326209418395, 0.9419777646489995, 0.04676784570357108);
\coordinate (B0) at (0.0, 0.40829947095934194, -0.08251603692661792);
\coordinate (C0) at (0.0, 1.3905867216880305, -0.26989735151234273);
\coordinate (D0) at (0.22699524986977337, 0.8159641094202493, -0.6138188687326387);
\fill[black!46!green!84, opacity=0.8] (A0) -- (B0) -- (C0) -- cycle;
\fill[black!26!green!47, opacity=0.9] (A0) -- (B0) -- (D0) -- cycle;
\fill[black!37!green!68, opacity=0.9] (B0) -- (C0) -- (D0) -- cycle;
\fill[black!6!green!11, opacity=0.9] (A0) -- (C0) -- (D0) -- cycle;
\draw [line width=0.05mm, color=black] (A0) -- (B0);
\draw [line width=0.05mm, color=black] (B0) -- (C0);
\draw [line width=0.05mm, color=black] (A0) -- (D0);
\draw [line width=0.05mm, color=black] (C0) -- (D0);
\draw [line width=0.05mm, color=black] (A0) -- (C0);
\draw [line width=0.05mm, color=black] (B0) -- (D0);

\coordinate (A1) at (1.078117894267925, 0.3049530456535753, 0.4964324257501687);
\coordinate (B1) at (0.6326146321737411, 0.7535620026926064, 0.1797672285342549);
\coordinate (C1) at (0.1871113700795573, 0.21988370900294885, 0.050483345904065884);
\coordinate (D1) at (0.8596098820435145, 0.17893939042482518, -0.16415428868604107);
\fill[black!46!green!84, opacity=0.8] (A1) -- (B1) -- (C1) -- cycle;
\fill[black!46!green!84, opacity=0.8] (A1) -- (C1) -- (D1) -- cycle;
\draw [line width=0.05mm, color=black, solid] (A1) -- (C1);
\fill[black!6!green!11, opacity=0.9] (A1) -- (B1) -- (D1) -- cycle;
\fill[black!15!green!27, opacity=0.9] (B1) -- (C1) -- (D1) -- cycle;
\draw [line width=0.05mm, color=black] (A1) -- (B1);
\draw [line width=0.05mm, color=black] (B1) -- (C1);
\draw [line width=0.05mm, color=black] (A1) -- (D1);
\draw [line width=0.05mm, color=black] (C1) -- (D1);
\draw [line width=0.05mm, color=black] (B1) -- (D1);

\coordinate (A2) at (0.539653053769138, 0.7987253449460457, 0.04295196090295757);
\coordinate (B2) at (0.09414979167495402, 0.26504705125638817, -0.08633192172723143);
\coordinate (C2) at (0.7666483036389113, 0.22410273267826453, -0.3009695563173384);
\coordinate (D2) at (0.3211450415447274, 0.6727116897172956, -0.6176347535332523);
\fill[black!45!green!81, opacity=0.8] (A2) -- (B2) -- (C2) -- cycle;
\fill[black!34!green!61, opacity=0.8] (A2) -- (B2) -- (D2) -- cycle;
\draw [line width=0.05mm, color=black, solid] (A2) -- (B2);
\fill[black!34!green!61, opacity=0.9] (B2) -- (C2) -- (D2) -- cycle;
\fill[black!6!green!11, opacity=0.9] (A2) -- (C2) -- (D2) -- cycle;
\draw [line width=0.05mm, color=black] (B2) -- (C2);
\draw [line width=0.05mm, color=black] (A2) -- (D2);
\draw [line width=0.05mm, color=black] (C2) -- (D2);
\draw [line width=0.05mm, color=black] (A2) -- (C2);
\draw [line width=0.05mm, color=black] (B2) -- (D2);

\coordinate (A3) at (0.09533800494530484, 0.1669579738068738, -0.22696307415914224);
\coordinate (B3) at (0.767836516909262, 0.1260136552287501, -0.44160070874924917);
\coordinate (C3) at (0.3223332548150782, 0.5746226122677811, -0.7582659059651631);
\coordinate (D3) at (0.5493285046848515, 0.0, -1.102187423185459);
\fill[black!26!green!47, opacity=0.8] (A3) -- (B3) -- (C3) -- cycle;
\fill[black!46!green!84, opacity=0.8] (A3) -- (B3) -- (D3) -- cycle;
\draw [line width=0.05mm, color=black, solid] (A3) -- (B3);
\fill[black!6!green!11, opacity=0.9] (B3) -- (C3) -- (D3) -- cycle;
\fill[black!37!green!68, opacity=0.9] (A3) -- (C3) -- (D3) -- cycle;
\draw [line width=0.05mm, color=black] (B3) -- (C3);
\draw [line width=0.05mm, color=black] (A3) -- (D3);
\draw [line width=0.05mm, color=black] (C3) -- (D3);
\draw [line width=0.05mm, color=black] (A3) -- (C3);
\draw [line width=0.05mm, color=black] (B3) -- (D3);

\end{tikzpicture}
  }
  \caption{Red-green refinement applied to a single macro element.
    From left to right:
    Unrefined element,
    regularly refined element (red),
    closure elements (green) to accommodate one, two and three new vertices, respectively.
  }
  \label{fig:rg3d}
\end{figure}

In 2\,d, each triangle is either refined regularly (\emph{red}) or bisected once (\emph{green}).
When using RG to produce a sequence of meshes, the green step must be reverted before applying the algorithm again, to prevent mesh degeneration.
In 3\,d, 
the 
procedure becomes significantly more complicated.
Following the refinement rules of Bey~\cite{bey},
nine different types of green refined elements occur, some of which have faces that are bisected twice.
We use a modified set of refinement rules, enforcing that all element faces also adhere to the rules of the 2\,d algorithm.
The procedure is summarized in \cref{alg:rg3d}.
Rather than nine types of irregular elements, this version of the algorithm yields only three different types, shown in \cref{fig:rg3d}.

\begin{algorithm}
  \footnotesize
  \begin{algorithmic}[1]
      \setstretch{1.5}
      \State Given a triangulation $\Tri$ and a set of elements marked for refinement $\Rc\subset\Tri$.
      \Repeat
        \State Apply regular (red) refinement to all $T\in\Rc$ \label{alg:rg:refCellsRed}
        \Repeat
          \State $\Rc_F\gets\emptyset$
          \ForAll{$T\in\Tri$}
            \State $\Rc_F\gets \Rc_F\cup\{{\footnotesize\text{faces }F\subset∂T\text{ with more than one hanging node on }∂F}\}$ \label{alg:rg:setupFacesRed}
          \EndFor
          \State Apply regular (red) refinement to all $F\in\Rc_F$  \label{alg:rg:refFacesRed}
        \Until{$\Rc_F = \emptyset$}
        \State $\Rc\gets\{{\footnotesize T\in\Tri\text{ with more than three hanging nodes on }∂T}\}$  \label{alg:rg:setupCellsRed}
      \Until{$\Rc = \emptyset$}
      \State $\Gc\gets\{{\footnotesize T\in\Tri\text{ with at least one hanging node on }∂T}\}$
      \State Apply irregular (green) refinement to all $T\in\Gc$
  \end{algorithmic}
  \caption{Red-Green Refinement (3\,d)}
  \label{alg:rg3d}
\end{algorithm}

Attempting to apply a refinement strategy like \cref{alg:rg3d} or other AMR techniques for unstructured grids
at scales of $\sim 10^{12}$ or more elements, poses a tremendous challenge.
In fact, even just storing the structural information needed for an unstructured mesh may be unfeasible.
We note here, that $10^{12}$ elements cannot be indexed with 32 bit integers.
iThus, if only a single 64 bit pointer were needed per such element, this alone adds 8 TByte to the global memory requirement.
In consequence, a typical data structure using multiple pointers to represent the topology of an
unstructured mesh and additionally storing the geometrical positions of the nodes, may lead to a prohibitive memory consumption.

Furthermore, any algorithm operating on data of such size will have to be of only linear (or almost linear) computational complexity,
excluding many traditional graph based algorithms from being used.
Even if the computational complexity is suitably low, the algorithm will still have to scale to many thousands of processor cores.
With this in mind, we argue that computations requiring such a high resolution and such a large problem size
are not feasible on conventional unstructured meshes combined with classical adaptive refinement strategies, since they
will become prohibitively expensive.
For such extreme scale problem sizes, using block-structured grids and exploiting their inherent advantages, including solvers of linear complexity, such as multigrid, will then become imperative.

Our FE-framework \hyteg is targeted primarily at these extreme scale problems that
can only be tackled using data structures and algorithms with minimal overhead.
For this, also the mesh information must be organized and stored efficiently and
such that the associated linear systems can be solved
with the best possible performance in a massively parallel setting.
In particular, as already explained in \cref{sec:hhg}, 
each process is responsible for the data corresponding to a 
subset $\Triloc\subset\Tri_0$, comprising only a few or even just one macro element.
Clearly, applying an AMR algorithm to this macro grid, which usually is several orders of magnitude smaller than the fine grid, is computationally relatively cheap, compared to solving the fine grid system.
Thus, a primary benefit of our construction is, that it suffices to run \cref{alg:rg3d} sequentially.
In fact, due to the small element count, parallelizing the algorithm would lead to extremely few loop iterations in lines \ref{alg:rg:refCellsRed} and \ref{alg:rg:refFacesRed}, while in lines \ref{alg:rg:setupFacesRed} and \ref{alg:rg:setupCellsRed}, refinement data from the neighboring elements is required.
The communication overhead would therefore be a major bottleneck.
The impact of using the refinement algorithm only in serial 
on the overall scalability will be discussed in \cref{sec:numex:scaling}.


  \section{A posteriori error estimate}\label{sec:error_estimate}
In the following sections an a posteriori error estimate for hierarchical meshes is derived.
The idea is based on the features and requirements on the specific linear solver suitable for HHG grids:
Since scalability is the central goal of our work, it is crucial to use solvers with optimal or near-optimal linear complexity, such as full multigrid solvers (FMG)~\cite{trottenberg2000multigrid}.
Computing a sequence of approximate and increasingly accurate solutions is therefore an inherent part of the solver.
In this setting, reusing these solutions of different accuracy to estimate the error seems to be a natural choice.
In contrast to other, more sophisticated methods to estimate the errors, this approach is not only extremely easy to integrate into our framework, it also comes at almost no additional cost.
In fact, since all the required data is already available, the workload for computing such an estimate almost exclusively comes from evaluating the norm, which can be reduced to a single vector dot product.

\subsection{Model problem and discretization}\label{sec:model_problem}
In order to derive a problem independent estimate, the problem description is kept abstract and general.
Consider an elliptic partial differential equation (PDE) over some domain $Ω\subsetℝ^d$, given in variational formulation:
Find $u\in V$ such that
\begin{equation}\label{eq:model_problem}
    a(u,v) = b(v) \quad\forall v\in V,
\end{equation}
where $V$ is a suitable Hilbert space of functions over $Ω$ (\eg, $V=H^1(Ω)$), $a\colon V⨯V→ℝ$ is a bilinear form, and $b\colon V→ℝ$ a linear form.

Let \cref{eq:model_problem} be discretized by conforming, nested finite element spaces $\SEQ{V_\ell}{\ell}$, \ie,
$V_{\ell-1} \subset V_{\ell} \subset V$ for all $\ell\in\SET{1,\dots,\L}$.
These spaces are associated with a sequence of triangulations $\SEQ{\Tri_\ell}{\ell}$ of $Ω$, as described in \cref{sec:hhg}.
Let $u∈V$ and $u_\ell\in V_\ell$ denote the solution of \cref{eq:model_problem} and corresponding discrete solution, and $e_\ell=(u-u_\ell)\in V$ the discretization error.
Furthermore, let $\|⋅\|$ be some norm over $V$.
Then, the grid convergence can be described by convergence factors
\begin{equation}\label{eq:conv_factor}
    \conv_\ell := \frac{\|e_{\ell}\|}{\|e_{\ell-1}\|}.
\end{equation}
Now suppose that, for given $\ell$, $\Tri_\ell$ resolves the problem sufficiently well to assert quasi-asymptotic convergence.
That means, there exist positive constants $ε\ll 1$ and $\conv$
such that
\begin{equation}\label{ass:asymptotic_convergence}
    \forall i\geq\ell\quad \conv \leq \conv_i \leq (1+ε)\conv < 1.
\end{equation}
Recursively applying \cref{eq:conv_factor} then allows estimating the fine grid error by
\begin{equation}\label{eq:err_est_anal}
\conv^{\L-\ell}\|e_\ell\| \leq \|e_{\L}\| \leq (1+ε)^{\L-\ell}\conv^{\L-\ell}\|e_\ell\|.
\end{equation}
While the above inequality is based on unknown quantities, it serves as a basis for a computable estimate ${η ≈ \|e_{\L}\|}$, that will be derived in the following section.

\begin{remark}\label{remark:discretization_order}
    $\conv$ is bounded from below by $2^{-q}$,
    where $q$ is the discretization order of the FE space with respect to $\|⋅\|$.
    With this optimal convergence factor, \cref{ass:asymptotic_convergence} is equivalent to the saturation assumption, see \eg~\cite{bank-saturation}.
\end{remark}

\subsection{Global error estimate}\label{sec:global_error_estimate}
We now introduce approximations $\errest{\ell}≈e_\ell$ and $\convest{\ell}≈\conv_\ell$, defined by
\begin{subequations}\label{eq:approximations}
\begin{equation}\label{eq:e_tilde}
    \errest{\ell} := e_\ell - e_\L = (u - u_\ell) - (u - u_\L) =  u_\L - u_\ell,
\end{equation}
\vspace{-1em}
\begin{equation}\label{eq:convest}
    \convest{\ell} := \frac{\|\errest{\ell}\|}{\|\errest{\ell-1}\|}.
\end{equation}
\end{subequations}
Let $0<j<\L$ and $\ell=\L-j$.
Then, replacing $e_\ell$ and $\conv_\ell$ in \cref{eq:err_est_anal} with these approximations leads to an a posteriori error estimate given by
\begin{equation}\label{eq:eta_glob}
    η_{j} := \convest{\ell}^{j}\|\errest{\ell}\|.
\end{equation}

\begin{theorem}\label{thm:eta_bounds}
    Let the error $e_\L$ be estimated by $η_j$ as defined in \cref{eq:eta_glob} and suppose that assumption \cref{ass:asymptotic_convergence} holds for $\ell=\L-j$.
    Then, $η_j$ is equivalent to the error norm, \ie, there exist positive constants $C_1$ and $C_2$ such that
    \begin{equation*}
        C_1 \|e_\L\| \leq η_j \leq C_2 \|e_\L\|.
    \end{equation*}
    Introducing the shorthand notation $\conv_ε := (1+ε)\conv$, these constants are given by
    \begin{equation*}
        C_{1} = (1+ε)^{-j}\frac{(1 - \conv_ε^{j})^{j+1}}{(1 + \conv_ε^{j+1})^{j}},\quad
        C_{2} = (1+ε)^{j}\frac{(1 + \conv_ε^{j})^{j+1}}{(1 - \conv_ε^{j+1})^{j}}.
    \end{equation*}
\end{theorem}

\begin{proof}
    Applying $\|⋅\|$ to \cref{eq:e_tilde},
    the triangle inequality implies
    \begin{equation*}
        \|e_\ell\| - \|e_\L\| \leq \|\errest{\ell}\| \leq \|e_\ell\| + \|e_\L\|.
    \end{equation*}
    Using \cref{eq:err_est_anal}, the contributions of $\|e_\L\|$ in the above inequality can be controlled by $\|e_\ell\|$ and vice versa, yielding
    \begin{align}
        \label{eq:e_tilde_bounds_ell}
        (1 - \conv_ε^{j})\|e_\ell\| \leq &\ \|\errest{\ell}\| \leq  (1 + \conv_ε^{j})\|e_\ell\|,\\
        \label{eq:e_tilde_bounds_L}
        (\conv_ε^{-j} - 1)\|e_\L\| \leq &\ \|\errest{\ell}\| \leq (\conv^{-j} + 1)\|e_\L\|.
    \end{align}
    Now, \cref{eq:e_tilde_bounds_ell} can be used to derive upper and lower bounds on the estimated convergence factor $\convest{\ell}$.
    In particular, by inserting \cref{eq:e_tilde_bounds_ell} into \cref{eq:convest}, we find
    \begin{equation*}
        \convest{\ell} = \frac{\|\errest{\ell}\|}{\|\errest{\ell-1}\|}
                        \leq \frac{\|e_{\ell}\|(1+\conv_ε^{j})}{\|e_{\ell-1}\|(1-\conv_ε^{j+1})}
                        \leq \conv_ε \frac{1+\conv_ε^{j}}{1-\conv_ε^{j+1}}.
    \end{equation*}
    Together with the analogously derived lower bound, one obtains
    \begin{equation}\label{eq:convest_bounds}
        \conv \frac{1-\conv_ε^{j}}{1+\conv_ε^{j+1}} \leq \convest{\ell} \leq \conv_ε \frac{1+\conv_ε^{j}}{1-\conv_ε^{j+1}}.
    \end{equation}
    Inserting \cref{eq:e_tilde_bounds_L} and \cref{eq:convest_bounds} into \cref{eq:eta_glob} yields
    \begin{equation*}
        \left((1+ε)^{-j} - \conv^j\right)\left(\frac{1-\conv_ε^{j}}{1+\conv_ε^{j+1}}\right)^j
        \leq \frac{η_j}{\|e_\L\|} \leq
        \left((1+ε)^{j} + \conv_ε^j\right)\left(\frac{1+\conv_ε^{j}}{1-\conv_ε^{j+1}}\right)^j.
    \end{equation*}
    $C_1$ is obtained by factoring out $(1+ε)^{-j}$ on the left-hand side.
    To keep the notation of $C_2$ consistent with $C_1$, the term $\left((1+ε)^{j} + \conv_ε^j\right)$ on the right-hand side is estimated from above by ${(1+ε)^{j}\left(1 + \conv_ε^j\right)}$.
\end{proof}

\begin{remark}\label{remark:preasymptotic}
    Provided $ε$ is sufficiently small,
    the accuracy of the estimate can be improved by increasing $j$.
    It can be easily verified that, for $ε=0$, both $C_1$ and $C_2$ approach one as $j\to ∞$.
    In practice though, increasing $j$ might often be problematic.
    Firstly, doing so implicitly makes $\Tri_\ell$ coarser, which makes it more challenging to control $ε$ in \cref{ass:asymptotic_convergence}.
    Secondly, for non-negligible $ε$, the constants actually deteriorate for increasing $j$ due to the factors $(1+ε)^{\pm j}$.
    Therefore, the choice of $j$ should be balanced to achieve the desired accuracy but retain the necessary robustness.
    This will be shown by an example in \cref{sec:numex:effectivity}.
\end{remark}
\begin{remark}\label{remark:algebraic_error}
    In the above analysis, only discretization errors are considered, \ie,
    it is assumed that the exact solutions $u_\ell$ and $u_\L$ are used in \cref{eq:e_tilde}.
    In practice though, these are usually perturbed by algebraic errors when the iterative solver is stopped and possibly
    also round-off errors which are of course also carried over to $\errest{\ell}$.
    For the above results to hold, it is hence necessary that these other error components do not exceed the discretization error.
    On the finest level $\L$, this should be satisfied anyway, as otherwise the accuracy of the discretization is not utilized.
    Here however, the requirement also extends to the lower levels $\ell$ and $\ell-1$. 
\end{remark}

\subsection{Local error indicator}\label{sec:error_indicator}
AMR is typically driven by some indicator for the local errors,
such that those elements can be marked for refinement that contribute most to the global error, see \eg~\cite{ainsworth1997}.
Defining the a local analogue to $η_j$ is straightforward:
For any subset $T\subsetΩ$, let $\norm{⋅}{T}$ be a norm on the local space $V\restrict{T}$,
such that $\norm{v}{T} = \|v\restrict{T}\|$ for all $v∈V$.
Local convergence factors $\conv_{\ell,T}$ and their approximations $\convest{\ell,T}$ are then defined as in \cref{eq:conv_factor,eq:convest}, respectively, but replacing $\|⋅\|$ with $\norm{⋅}{T}$.
Then, an estimate to $\norm{e_\ell}{T}$ is given by
\begin{equation}\label{eq:eta_loc_j}
    η_{j,T} := \convest{\ell,T}^j\norm{\errest{\ell}}{T}.
\end{equation}
Assumption \cref{ass:asymptotic_convergence} must also be adapted accordingly:
Suppose there exist positive constants $ε\ll 1$ and $\conv_{T}$ for all $T\in\Tri_0$, such that
\begin{equation}\label{ass:local_asymptotic}
    \forall i\geq\ell\quad \conv_T \leq \conv_{i,T} \leq (1+ε)\conv_T < 1.
\end{equation}
It is important to note that $\conv_T$ must be expected to vary significantly with respect to $T$ -- otherwise there would be no need for AMR in the first place.
Global upper and lower bounds on $\conv_{i,T}$ are given by
\begin{equation*}\label{eq:min_max_conv}
    \convmin := \min_{T∈\Tri_0} \conv_T,\quad
    \convmax := (1+ε)\max_{T∈\Tri_0} \conv_T.
\end{equation*}

Since \cref{eq:eta_loc_j} is canonical with \cref{eq:eta_glob}, most of the analysis in \cref{sec:global_error_estimate} can be reused.
However, one has to keep in mind, that the objective of a local estimate is quite different from a global one.
Global effectivity indices should be close to one, \ie, $γ=\sfrac{η}{\|e_\L\| ≈ 1}$.
For the local indices on the other hand, it is sufficient if there is a constant $C$ such that ${γ_T=\sfrac{η_{T}}{\norm{e_\L}{T}}≈C}$ for all $T∈\Tri_0$,
or, equivalently,
\begin{equation*}\label{eq:epsilon_loc_opt}
    ∂γ := \frac{\max_{T∈\Tri_0} γ_T}{\min_{T∈\Tri_0} γ_T} ≈ 1.
\end{equation*}

Suppose the local estimates are equivalent to the error, \ie,
for each element $T\in\Tri_0$ there exist positive constants $C_{1,T}$, $C_{2,T}$ such that
\begin{equation}\label{eq:bounds_local_effectivity}
    C_{1,T}\norm{e_\L}{T} \leq η_{T} \leq C_{2,T}\norm{e_\L}{T}.
\end{equation}
Then, it is easy to see that an upper bound for $∂γ$ is given by
\begin{equation}\label{eq:delta_bound_c12}
    ∂γ \leq \frac
        {\max_T C_{2,T}}
        {\min_T C_{1,T}}.
\end{equation}
For the proposed local error estimator, these bounds are trivially given by the constants $C_1$ and $C_2$ from \cref{thm:eta_bounds}.

\begin{corollary}\label{thm:bounds_eta_jT}
    Let the local errors be estimated by $η_{j,T}$ as defined in \cref{eq:eta_loc_j} and suppose assumption \cref{ass:local_asymptotic} holds for $\ell=\L-j$.
    Then, $∂γ$ is bounded by
    \begin{equation*}
        ∂γ \leq
        (1+ε)^{2j}
        \left(\frac{1 + \convmax^{j+1}} {1 - \convmax^{j+1}}\right)^{j}
        \left(\frac{1 + \convmax^{j}}   {1 - \convmax^{j}}  \right)^{j+1}.
    \end{equation*}
\end{corollary}
\begin{proof}
    As $η_{j,T}$ is defined analogous to $η_j$, \cref{thm:eta_bounds} can be employed element-wise to derive upper and lower bounds as in \cref{eq:bounds_local_effectivity}.
    In particular, $C_{1,T}$ and $C_{2,T}$ in \cref{eq:delta_bound_c12} are obtained by replacing $\conv$ in $C_1$ and $C_2$ with $\conv_T$.
    Since the constants converge monotonically as $\conv_T$ decreases, both ${\min_T C_{1,T}}$ and ${\max_T C_{2,T}}$ are defined by the worst local convergence factor and it holds
    \begin{equation*}
        ∂γ \leq \frac
            {(1+ε)^j\sfrac{(1 + \convmax^{j})^{j+1}}{(1 - \convmax^{j+1})^{j}} }
            {(1+ε)^{-j}\sfrac{(1 - \convmax^{j})^{j+1}}{(1 + \convmax^{j+1})^{j}} }.
    \end{equation*}
    Reordering the terms yields the 
    result.
\end{proof}
\begin{remark}\label{remark:eta_jT}
    As for the global estimate, the accuracy of this local estimate can be improved by increasing $j$.
    However, the same considerations as in \cref{remark:preasymptotic} apply.
\end{remark}

The above result shows that the reliability of $η_{j,T}$ heavily depends on $\convmax$.
In particular when considering problems with 
singularities, as caused, e.g.~by a reentrant corner, the local convergence order close to the singularity may be very poor.
Adding to this, one might be forced to keep $j$ small (s. \cref{remark:preasymptotic}), further increasing the upper bound on $∂γ$.
To better deal with these kind of issues, we propose to omit the scaling factor in \cref{eq:eta_loc_j}, yielding
\begin{equation}\label{eq:eta_loc_prime_j}
    η'_{j,T} := \norm{\errest{\L-j}}{T}.
\end{equation}

\begin{corollary}\label{thm:bounds_eta_prime_jT}
    Let the local errors be estimated by $η'_{j,T}$ as defined in \cref{eq:eta_loc_prime_j} and suppose assumption \cref{ass:local_asymptotic} holds for $\ell=\L-j$.
    Then, $∂γ$ is bounded by
    \begin{equation*}
        ∂γ \leq \frac
        {\convmin^{-j}+1 }
        {\convmax^{-j}-1 }
    \end{equation*}
\end{corollary}
\begin{proof}
    Here, the constants $C_{1,T}$, $C_{2,T}$ can be taken from \cref{eq:e_tilde_bounds_L},
    again, replacing $\conv$ by $\conv_T$.
    Note that, as above, $C_{1,T}$ is minimized for the element with the worst local convergence.
    The maximum of $C_{2,T}$ on the other hand is achieved for the best local convergence factor.
\end{proof}
\begin{remark}\label{remark_eta_prime_jT}
    While $η'_{j,T}$ is less sensitive with respect to large $\convmax$, than $η_{j,T}$,
    the upper bound on $∂γ$ here also depends on the best local convergence factor $\convmin$.
    In particular, the accuracy depends on the ratio between worst and best convergence factors $\sfrac{\convmax}{\convmin}$.
    Considering \cref{remark:discretization_order}, this implies that $η'_{j,T}$ is more suitable for lower order FE.
    Also noteworthy,
    instead of converging as $j\to∞$, here $∂γ$ diverges at a rate of $\bigO((\sfrac{\convmax}{\convmin})^j)$.
\end{remark}

The implications of \cref{thm:bounds_eta_jT,thm:bounds_eta_prime_jT} on the choice between $η_{j,T}$ and $η'_{j,T}$, as well as the value of $j$ are summarized below.
The major factor is of course the discretization order $q$ of the chosen FE space.
\paragraph{High order} For $q>3$, the scaled estimator $η_{j,T}$ is usually preferable.
    Although $∂γ$ can then be improved by increasing $j$, it may still be beneficial to choose $j=1$.
    This is due to the fact that convergence is expected to reach asymptotic behavior later when using higher order methods.
\paragraph{Moderate order}\label{par:moderate_order} In case of $q=2,3$, both estimators are similarly suitable.
    While $η_{j,T}$ is more accurate if $\convmax≈\convmin$, the advantage shifts towards $η'_{j,T}$ as $\convmax$ increases.
    In particular if $\convmax > \sfrac{1}{2}$, \ie, there are regions with sublinear convergence, the unscaled estimator is clearly preferable.
    Another major advantage of $η'_{j,T}$ is that is most accurate for $j=1$, while $η_{j,T}$ relies on less robust choices of $j$ to be competitive.
\paragraph{Low order} For $q=1$, the unscaled estimator $η'_{j,T}$ clearly outperforms $η_{j,T}$ in every regard.

\subsection{Implementation in \hyteg}\label{sec:errest_impl}
An efficient implementation the above defined error estimator in \hyteg is straightforward.
Since the framework is built for multigrid methods, the computational cost can be minimized by reusing the coarse grid solutions, that are an inherent by-product of applying an FMG solver (see \eg \cite{trottenberg2000multigrid}).

As usual, the discretization of \cref{eq:model_problem} on level $\ell$, is formulated as linear system
\begin{equation}\label{eq:lse}
    A_\ell \uh_\ell = \bh_\ell
\end{equation}
with $n_\ell = \dim V_\ell$ unknowns.
When using an FMG to solve the system on level $\L$, solutions $\uh_\ell$ for all levels $\ell<\L$ are also available at no extra cost.
Computing $\errest{\ell}$ then only requires interpolating the coarse solutions to the finest grid level.
Using the isomorphism between $V_\L$ and $ℝ^{n_\L}$, the estimated error can be formulated as
\begin{equation*}\label{eq:errest_Rn}
    V_\L \ni \errest{\ell} = u_\L - u_\ell ≡ \uh_\L - I_{\ell}^{\L}\uh_\ell =: \eest_\ell \in ℝ^{n_\L},
\end{equation*}
where $I_{\ell}^{\L}$ is an appropriate interpolation operator.
The operation count associated with the above equation can be further reduced by setting $\ell=\L-1$ and using the same interpolation operator as in the multigrid scheme.
Then, applying $I_{\ell}^{\L}$ is essentially free since, for all $\ell$, the solution $\uh_\ell$ is already being interpolated to $\ell+1$ as part of the FMG.
The only extra cost is storing the interpolated values in a separate data structure, as they would otherwise be overwritten in the next step of the algorithm (see \cref{alg:fmg}).

Regarding the norm to evaluate $η_j$,
the most obvious choices are $L^2$ and energy norm.
In the experiments presented in \cref{sec:numex}, the $L^2$-norm is used, \ie, ${\|⋅\| = \|⋅\|_{L^2(Ω)}}$ and ${\|⋅\|_T = \|⋅\|_{L^2(T)}}$.
Sine the estimate lives in the discrete space $V_\L$,
computing the norm simplifies to
$\|\errest{\ell}\|^2 = \eest_{\ell}\t M_\L\eest_{\ell}$,
where $M_\L$ is the standard mass matrix.
Using only standard operators, automatically generated and highly optimized kernels can be utilized to compute $η_j$.

It is worth mentioning that, in \hyteg, each macro-vertex, -edge, -face, and -cell owns the data corresponding to the DoF in its interior.
The inter-macro boundaries are treated as ghost layers.
Matrix-vector products are then computed piecewise on the volume macros only (faces in 2\,d, cells in 3\,d), which requires no communication.
The data on the interfaces is obtained by additively collecting the data from the ghost layers of adjacent volume elements.
This communication step is not required here, since the local $L^2$-norms for $η_T$ only require the local data on the ghost nodes.
The global estimate can then be obtained as the sum of squared local estimates.
The entire procedure and how it is interleaved in the FMG solver is summarized in \cref{alg:fmg}.

As of now, the implementation has only been realized for linear Lagrange FE, where
optimal $L^2$-convergence is quadratic ($q=2$).
Following the guidelines in \cref{sec:error_indicator}, the unscaled estimates are employed, \ie,
    $η_{T} := η'_{1,T} = \norm{\errest{\L-1}}{T}$.

\begin{remark}\label{remark:norm}
    Switching to a different norm is straightforward, \eg,
    $M_\L$ can be replaced with the stiffness matrix $A_\L$ to obtain the energy norm.
    The mass matrix has the advantage of being piecewise constant, which leads to extremely efficient matrix-free kernels.
    Mass lumping~\cite{hughes1978mass} could be used to further reduce the computational cost.
    Regarding the local estimates, one could think of even more advanced approaches, such as goal oriented error estimates, see \eg~\cite{becker2001optimal,richter2015goal}.
    Suppose the goal function can be described by a piecewise linear operator $G:V→W$, where $W$ is a vector space with norm $\norm{⋅}{W}$.
    Then, the analysis in \cref{sec:error_indicator} holds for a local (semi)-norm, defined by $\norm{v}{T} = \norm{G(v\restrict{T})}{W}$.
\end{remark}

\begin{algorithm}
    \footnotesize
    \begin{algorithmic}[1]
        \setstretch{1.5}
        \State{Given $\SEQ{A_\ell, \bh_\ell}{\ell=0}^{\L}$, zero initialized vectors $\SEQ{\uh_\ell}{\ell=0}^{\L}$, $\SEQ{\wh_\ell}{\ell=0}^{\L-1}$ and $ν∈ℕ$.}
        \Procedure{FMG}{$\SEQ{A_\ell, \uh_\ell, \bh_\ell}{\ell},\ \SEQ{\wh_\ell}{\ell},\ ν$}
            \For{$\ell=0,\dots,\L$}
                \For{$i=1,\dots,ν$}
                    \State{$\uh_\ell \gets \textproc{V-cycle}(\SEQ{A_j, \uh_j, \bh_j}{j=0}^{\ell})$}
                \EndFor
                \Comment{{\footnotesize Solve for $\uh_\ell$ using $ν$ MG V-cycles}}
                \If{$\ell < \L$}
                    \State{$\uh_{\ell+1} \gets I_{\ell}^{\ell+1}\uh_\ell$}
                    \Comment{{\footnotesize Prolongate $\uh_\ell$ to obtain initial guess on $\ell+1$}}
                    \State{$\wh_{\ell} \gets \uh_{\ell+1}$}
                    \Comment{{\footnotesize Store prolongated coarse grid solution $\uh_\ell$}}
                \EndIf
            \EndFor
        \EndProcedure

        \Procedure{ErrorEstimate}{$\uh_\L,\ \SEQ{\wh_\ell}{\ell},\ j$}\label{alg:proc:errest}
            \State{$\ell=\L-j$}

            \State{$\wh_{\ell-1} \gets I_{\ell}^{\L}\wh_{\ell-1}$}
            \State{$\wh_{\ell} \gets I_{\ell+1}^{\L}\wh_{\ell}$}
            \Comment{{\footnotesize Prolongate coarse grid solutions to fine grid}}

            \State{$\|\errest{\ell-1}\| \gets \sqrt{(\uh_{\L} - \wh_{\ell-1})\t M_\L (\uh_{\L} - \wh_{\ell-1})}$}
            \Comment{{\footnotesize Coarse grid error estimates}}
            \State{$\|\errest{\ell}\| \gets \sqrt{(\uh_{\L} - \wh_{\ell})\t M_\L (\uh_{\L} - \wh_{\ell})}$}
            \Comment{{\footnotesize ($η_T$ obtained as by-product of $\|\errest{\L-1}\|$)}}


            \State{$\convest{\ell} \gets \sfrac{\|\errest{\ell}\|}{\|\errest{\ell-1}\|}$}
            \Comment{{\footnotesize Convergence estimate}}
            \State{$η_{j} \gets \convest{\ell}^{j} \|\errest{\ell}\|$}
            \Comment{{\footnotesize Fine grid error estimate}}
        \EndProcedure
    \end{algorithmic}
    \caption{FMG solver and error estimator}
    \label{alg:fmg}
\end{algorithm}
  \section{Numerical experiments}\label{sec:numex}
In the following, the concepts introduced in the previous sections will be applied to some 
illustrative examples.
The main goal is to investigate the effectivity of the local and global error estimates as well as the convergence behavior achieved with the resulting meshes.
To that end, Poisson problems with different requirements on mesh adaptivity are considered, \ie, \cref{eq:model_problem} is defined by
\begin{equation}\label{eq:poisson}
    a(u,v) = ∫_Ω ∇v⋅∇u,\quad
    b(v) = ∫_Ω v f\quad\text{for some given } f.
\end{equation}
As stated in \cref{sec:error_estimate}, our results are not limited to this class of problems.
However, for Poisson's equation it is 
easy to define \emph{manufactured} solutions with well-defined requirements on mesh adaptivity.
This is essential for a thorough 
numerical evaluation of our method.

The discrete system \cref{eq:lse} is solved using an FMG method with
Gauss Seidel iterations for smoothing and a conjugate gradient (CG) method to solve the the coarse grid system ($\ell=0$).
In \cref{sec:numex:waves,sec:numex:LShape}, the focus is on the convergence of \cref{alg:kplusl,alg:kl}.
Therefore, in order to keep the results free from any influences of the error estimator, the refinement process will be driven by the actual local $L^2$ errors, rather than $η_T$.

The estimator will then be examined in \cref{sec:numex:effectivity}.
In particular, this includes a study of its robustness with respect to the parameters $j$ and $ν$.
Note that the coarse solutions $u_{\ell}$, required to compute $\errest{\ell}$, are not computed exactly but only approximately by applying a few V-cycles within the FMG algorithm.
Clearly, the associated algebraic errors may have an impact on the effectivity of $η_j$ and $η_T$ (s.~\cref{remark:algebraic_error}).
On the other hand, using more V-cycles within FMG increases the computational cost.

In all of the following experiments, to ensure that 
the algebraic errors on the finest level are sufficiently small, additional V-cycles are applied on level $\L$, until the error $\|u-u_\L\|$ has converged, \ie, the first two significant digits do not change any more.
Note that this is 
computed after FMG terminates, and that the coarse solutions $u_{\L-j}$ are not updated during this process.


\subsection{Convergence study for a simple problem requiring AMR}\label{sec:numex:waves}
There are various different reasons why mesh adaptivity may be required for a given problem.
In this section, we consider a problem with a smooth but highly oscillating solution, as described in \cref{ex:waves}.

\begin{example}\label{ex:waves}
    Let the domain be given by the unit square $Ω = (0,1)^2 \subset ℝ^2$, and let the right hand side $f$ of \cref{eq:poisson} be defined such that the analytic solution satisfies
    \begin{equation}\label{eq:waves_solution}
        u(x,y) = w(x)w(y),
        \quad\text{where}\quad
        w(t) = 1 - \cos\left(e^{α(t-1)} ω t\right).
    \end{equation}
    Note that, albeit being smooth, the problem requires a significantly higher resolution in some parts of the domain than in others, since the frequency of the oscillations in $u$ increases exponentially towards the upper right corner of $Ω$.
    Thus also the gradients and higher derivatives of the true solution are higher in this region, but remain bounded.
    The wave number is given by $\sfrac{ω}{2π}$, while $α$ controls how much the oscillations are concentrated in the corner.
    Consequently, the parameters $ω$ and $α$ can be used to control the requirement on the global resolution and the demand for adaptivity, respectively.
    As a reference, \cref{fig:waves2d_u} shows a plot of $u$ over $Ω$ with moderate parameters $α$ and $ω$.
    In the experiments, more challenging values ${ω=16π},\ {α=10}$ are being used.
    \begin{figure}
        \centering
        \includegraphics[width=0.7\textwidth, trim=0 800 0 950, clip]{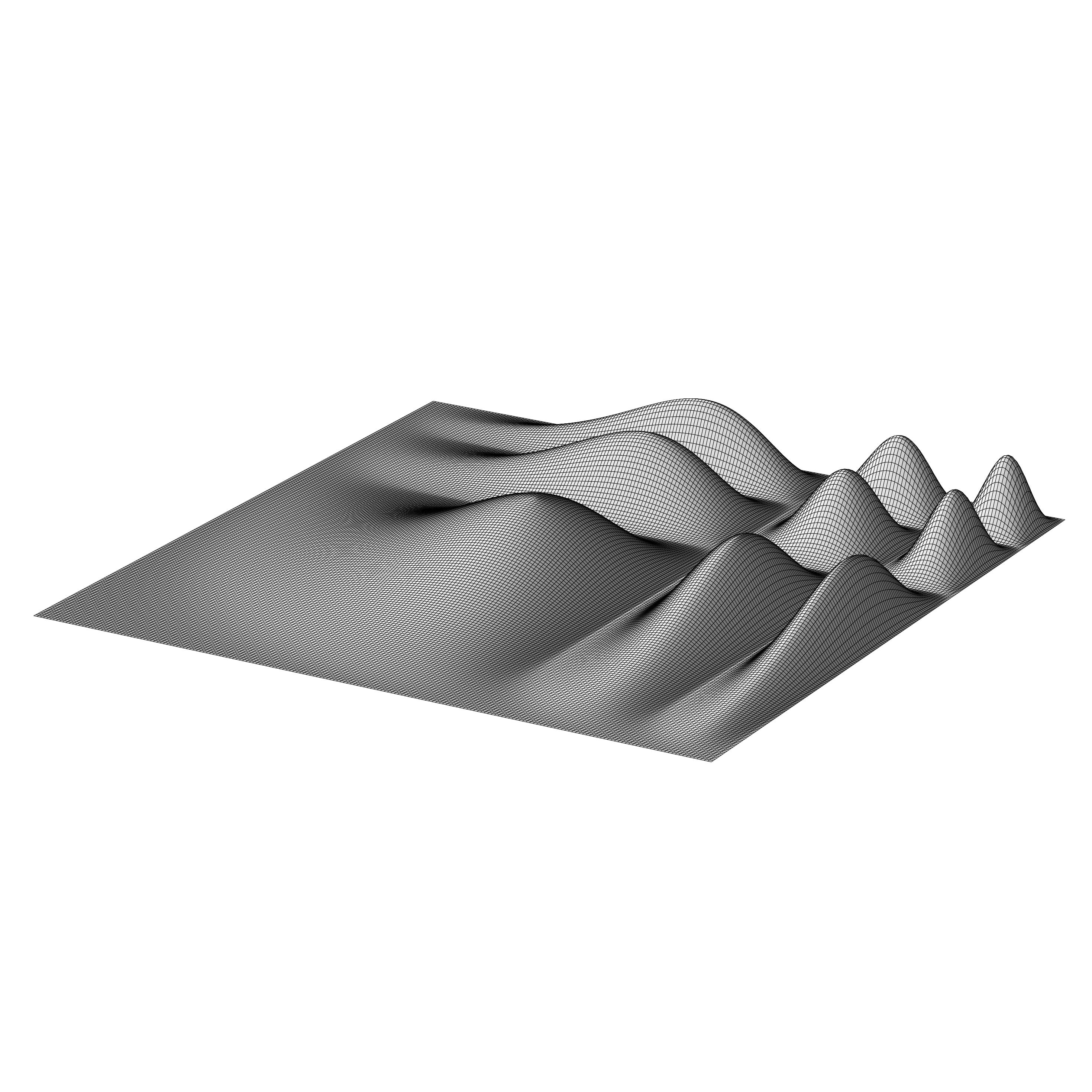}
        \caption{
            Solution $u$ as defined by \cref{eq:waves_solution} with $α=2$ and $ω=6π$.
        }
        \label{fig:waves2d_u}
    \end{figure}
\end{example}

An optimal sequence of grids should result in an error convergence as given by the discretization order.
For first order FE and the $L^2$-norm, that means $\bigO(\hmean^2)$, where $\hmean$ is a quantity such that $\hmean^d$ is proportional to the number of DoF in a given mesh $\Tri$.
We therefore define $\hmean := \left(\sfrac{|Ω|}{\#\Tri}\right)^{\sfrac{1}{d}}$.
Now, let $\hmean$ and $\Hmean$ denote the above defined mean mesh width of $\Tri_\ell^k$ and $\Tri_0^0$, respectively.
Furthermore, let $e_{\hmean}$ and $e_{\Hmean}$ denote the errors of the respective grids.
Then, for sufficiently small $\hmean$, the ratio
\begin{equation}\label{eq:errratio}
    \frac{\Hmean^2\norm{e_{\hmean}}{}} {\hmean^2\norm{e_{\Hmean}}{}} =: \errratio
\end{equation}
should be approximately constant with respect to $\hmean$.
The value $\errratio$ also poses a metric for the quality of a mesh sequence: The smaller $\errratio$ is, the more accuracy improvement per DoF.

\begin{figure}
    \centering
    \input{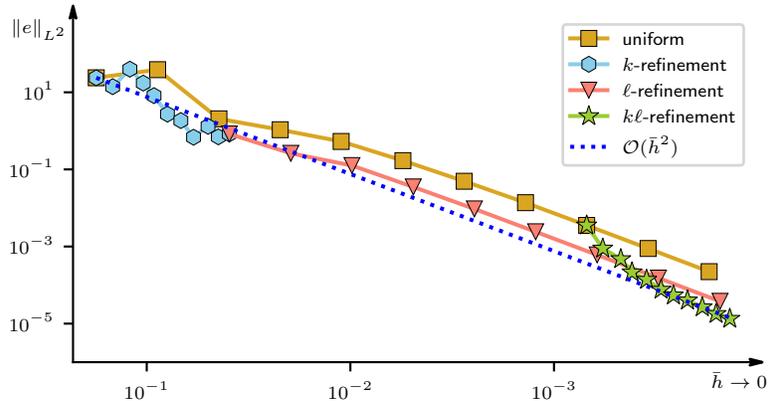}
    \caption{
        $L^2$ errors for different refinement schemes applied to \cref{ex:waves}.
        The yellow squares correspond to uniform refinement, the blue hexagons and red triangles to the loops over $k$ and $\ell$ in \cref{alg:kplusl}, respectively, and the green stars to the loop over $k$ in \cref{alg:kl}.
        The dotted blue graph marks the desired quasi optimal behavior with $\errratio=1$.
    }\label{fig:waves2d_l2err}
\end{figure}

The initial coarse grid $\Tri^0_0$ is a uniform triangulation comprising 32 macro elements.
As a baseline, the problem is solved using only uniform refinement.
The resulting convergence is displayed in \cref{fig:waves2d_l2err}.
Asymptotically, the optimal convergence of $\mathcal{O}(h^2)$ is achieved, but only once the oscillations have been resolved sufficiently well.
In the beginning of the refinement process, we observe a pre-asymptotic behavior with significantly slower convergence, yielding a rather large value of $\errratio=9.8$ in the finest grid.

This can be alleviated by using a coarse grid that captures the features of $u$ better.
In particular, regions with oscillations of higher frequency require a higher mesh resolution.
Note that, if the oscillations are resolved sufficiently well by the coarse grid,
the fine grid error can be expected to be reduced at the asymptotic rate.
With such a coarse grid, HHG meshes are expected to perform very well for this problem.
Such a grid can be obtained by either \cref{alg:kplusl}, or \cref{alg:kl}.
For both approaches, we use $\K=10$ adaptive refinement steps, where in each step $k$, $10\%$ of the macro elements in $\Tri^{k}$ are marked for \emph{red} refinement.
The fine grids are constructed as described in \cref{sec:hhg} with finest grid level $\L=8$.
We point out that the example is chosen such that if the coarse mesh $\Tri^\K_0$ has a suitable balance of the errors,
then this is expected to carry over to the final mesh $\Tri^\K_\L$.
In other words, here, once a suitable $\Tri^\K_0$ has been found, then uniform refinement leads to optimal convergence.

As expected, the final grids $\Tri^\K_\L$ of both algorithms yield significantly better results than uniform refinement,
with $\errratio=1.5$ and $\errratio=0.9$, respectively.
Since the improvement for \cref{alg:kl} is quite a bit better than for \cref{alg:kplusl}, it becomes clear that
local fine grid errors $\|e_\L\|_T$ are indeed a more effective refinement criterion than coarse grid errors $\|e_0\|_T$.
A convergence analysis is also shown in \cref{fig:waves2d_l2err}.

\subsection{Convergence for a problem with singularity}\label{sec:numex:LShape}
While the previous example shows that AMR can be beneficial even in case of smooth problems, this is even more so
when dealing with singularities, such as the corner singularity in \cref{ex:lshape}.

\begin{example}\label{ex:lshape}
    Let $(r, φ)$ denote polar coordinates corresponding to $\xx ∈ ℝ^2$, and let the domain and analytic solution be given by
    \begin{equation*}
        Ω=\SET{\xx∈(-1,1)^2 \setsep φ∈\left(0,\tfrac{3}{2}π\right)}
        \quad\text{and}\quad
        u(r,φ) = r^{\sfrac{2}{3}} \sin\left(\tfrac{2}{3}φ\right).
    \end{equation*}

    \begin{figure}
        \centering
        \subfloat[\label{fig:lshape_mesh_kPlusL}\cref{alg:kplusl}]
        {
             \includegraphics{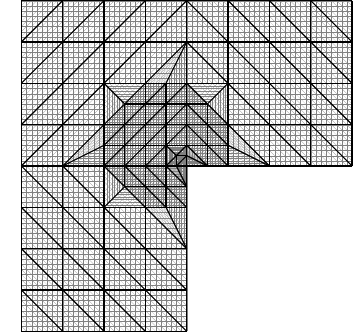}
        }
        \hfill
        \subfloat[\label{fig:lshape_mesh_kl}\cref{alg:kl}]
        {
             \includegraphics{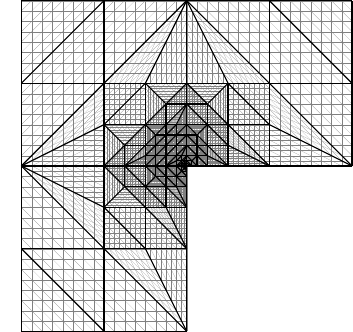}
        }
        \caption{
            Triangulations $\Tri_\ell^k$ of the L-shape domain, produced by different refinement procedures.
            Note that $k=5$ and $\ell=3$ is chosen for visualization purposes only and this particular combination does not actually come up in either of the two schemes:
            \cref{alg:kplusl} provides meshes $\Tri^0_0, \dots, \Tri^{10}_0, \dots, \Tri{10}^8$,
            while \cref{alg:kl} only produces $\Tri_8^{0}, \dots, \Tri_8^{10}$.
            Since $\Tri_8^{10}$, the configuration that comes up in both schemes, is not suitable for printing, $\Tri_3^5$ is shown instead.
            }
            \label{fig:lshape_mesh}
        \end{figure}
\end{example}

It is well known that uniform refinement of this \emph{L-shape} problem leads to $L^2$-convergence in $\bigO(h^{\sfrac{4}{3}})$.
This is true even when the boundary conditions and right hand side of the problem are smooth functions.
There are various techniques to restore optimal convergence order, such as modifying
the approximation in a neighborhood of the singularity, \eg~\cite{babuska1972, zenger1978, egger2014energy},
or enriching the FE space with non-smooth basis functions, see \eg~\cite{strang2008,grisvard2011}.
A more commonly used approach is the use of gradually refined meshes $\SEQ{\Tri^k}{k}$, instead of uniform refinement.
According to \cite{gradedMR}, these grids must satisfy the following condition:
Let $h_T$ and $r_T$ denote the mesh width of an element $T$ and its distance from from the singularity, respectively,
and let $h_{\max}=\max_{T\in\Tri^k} h_T$.
For all $T∈\Tri_k$, it must hold
\begin{equation}\label{eq:lshape_grading}
    h_{\max} r_T^{\sfrac{1}{3}} \lesssim h_T \lesssim  h_{\max} r_T^{\sfrac{1}{3}}.
\end{equation}

Clearly, satisfying such a condition on a block-structured grid is not easily possible, as the following experiment shows.
Analogous to the previous section, the domain is discretized by a very coarse, uniform initial grid $\Tri^0_0$, comprising only 24 triangles.
Again, this grid is used as basis for both \cref{alg:kplusl} and \cref{alg:kl}.
The parameters $\K$, $\L$, and refinement criterion are identical to the experiment in \cref{sec:numex:waves}.

\begin{figure}
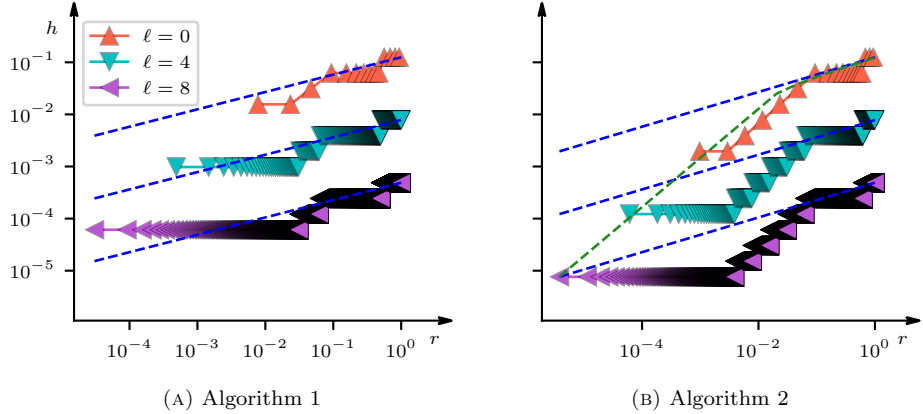

    \centering
    \subfloat[\label{fig:lshape_refinement_kPlusL}\cref{alg:kplusl}]
    {\input{figures/Lshape_refinement_comparison.pgf}}
    \subfloat[\label{fig:lshape_refinement_kl}\cref{alg:kl}]
    {\input{figures/Lshape_refinement.pgf}}
    \caption{
        Mesh width $h$ with respect to the distance $r$ for the meshes $\Tri_\ell^\K$ from the above experiment on \cref{ex:lshape}.
        The blue and green dashed lines mark the mesh requirements given in \cref{eq:lshape_grading,eq:lshape_grading_kl}, respectively.
    }
    \label{fig:lshape_refinement}
\end{figure}

The result for \cref{alg:kplusl} shown in \cref{fig:lshape_l2err} is as expected.
While the adaptive $k$-refinement yields the desired $\bigO(\hmean^2)$ convergence, the order reduces
to $\bigO(\hmean^{\sfrac{4}{3}})$ when using a refined coarse grid $\Tri^\K_0$ and continuing with uniform refinement.
This is no surprise, since the objective of the adaptive refinement procedure was minimizing the error for the coarse grid,
not the associated HHG constructed by uniform refinement on top of it.
The impact of this becomes even more obvious when analyzing the mesh parameters $h_T$ and $r_T$ of the fine grids.
While $\Tri^\K_{\ell}$ adheres quite well to \cref{eq:lshape_grading} for $\ell=0$, the same is not true at all for $\ell>0$.
In particular, these meshes are too coarse near the singularity, as can be seen in~\cref{fig:lshape_refinement_kPlusL}.

\begin{figure}
    \centering
    \subfloat[\label{fig:lshape_l2err}\cref{alg:kplusl}]
    {\input{figures/Lshape_l2err_err_V6_kpl.pgf}}
    \subfloat[\label{fig:lshape_kl_l2err}\cref{alg:kl}]
    {\input{figures/Lshape_l2err_err_V6_kl_klsep.pgf}}
    \caption{
        $L^2$ errors for different refinement schemes applied to \cref{ex:lshape}, analogous to \cref{fig:waves2d_l2err}.
        In addition to the desired $\bigO(\hmean^2)$ convergence, the expected convergence rate for uniform refinement of $\bigO(\hmean^{\sfrac{4}{3}})$ is marked by the orange dotted graph.
        To better compare the results in \subref{fig:lshape_kl_l2err} to those in \subref{fig:lshape_l2err}, the convergence rates with respect to $k$ and $\ell$ in \cref{alg:kl} are also shown individually.
        The gray graphs are obtained by solving the problem on $\Tri^k_0$ for all $k$, and $\Tri^\K_\ell$ for all $\ell$, using the coarse grids $\Tri^k$ arising in \cref{alg:kl}.
    }
\end{figure}

Since the number of uniform refinement steps $\L$ is known in advance,
a requirement on the coarse grid $\Tri^k$ can be derived, such that that the upper bound in \cref{eq:lshape_grading} is satisfied on the fine grid.
Let $\ell\leq\L$, $T_\ell\in\Tri^k_\ell$, and $T\in\Tri^k_0$ such that $T_\ell\subset T$.
Furthermore, let $h_{T,\ell}$ and $r_{T,\ell}$ denote the respective parameters of $T_\ell$.
It is easy to verify that
$h_{T,\ell} = 2^{-\ell} h_T$, and $r_{T,\ell} \geq r_T - \tfrac{1}{2} h_T + \tfrac{1}{2} h_{T,\ell}$.
Consequently, $\Tri_\ell^k$ satisfies the upper bound in \cref{eq:lshape_grading}, if, for all $T∈\Tri_0^k$, it holds
\begin{equation}\label{eq:lshape_grading_kl}
        h_T \lesssim h_{\max} \left(r_T - \frac{1-2^{-\L}}{2}h_T\right)^{\sfrac{1}{3}}.
\end{equation}
While 
interesting, this relation less useful in practice, since it is specific to this particular problem.
However, the objective of this article is to provide a fully automatic refinement scheme
that is applicable to general problems, where the requirements on the mesh are not known a priori.

As it turns out, artificially constructing meshes adhering to \cref{eq:lshape_grading_kl} is not even necessary, since $k\ell$-refinement provides just the right kind of mesh sequence.
\Cref{fig:lshape_refinement_kl} shows that the resulting coarse grids indeed satisfy \cref{eq:lshape_grading_kl}.
Moreover, the fine grids $\Tri^\K_\ell$ all adhere to the upper bound in \cref{eq:lshape_grading}.
Of course, for larger $\ell$, these graphs 
will increasingly take the form of staircases.
Consequently, the lower bound in \cref{eq:lshape_grading} is violated even more severely than before.
However, this has only a minor effect on the grid convergence, as can be seen in \cref{fig:lshape_kl_l2err}.
Appropriate post processing reveals that the convergence with respect to $\ell$ is in fact very close to the optimal rate, with $\errratio=1.63$ on $\Tri^\K_\L$.
In 2\,d, this means that only 1.6 times more DoF are required to achieve the same accuracy as an optimal mesh would.
Such a moderate factor should be easily compensated and outweighed by the superior performance of HHG
compared to a fully unstructured grid.
As a reference, we point to the effectiveness of the optimizations applied in \hyteg~\cite{boehm2024}, which are typically resulting in a speedup of an order of magnitude or more.
Many of these optimizations rely on and are only possible for hierarchically structured meshes of HHG type.
While these results have been found for conventional CPU types, one may expect that the performance
advantages of HHG grids will become even more pronounced on modern GPU-based supercomputers.
To give a more tangible visualization of the difference between the grids produced by \cref{alg:kplusl,alg:kl}, some examples are shown in \cref{fig:lshape_mesh}.

\subsection{Effectivity of the error estimate}\label{sec:numex:effectivity}
%
%
A widely used metric for error estimators is the so called effectivity index, defined as
${\sfrac{η_j}{\|e_\L\|}}$, which ideally should be close to one.
The analysis in \cref{sec:error_estimate} gives upper and lower bounds on the effectivity indices of our approach, namely the constants $C_1$ and $C_2$.
In order to validate these results, error estimates arising in the experiments from \cref{sec:numex:waves,sec:numex:LShape} are investigated for varying parameters $j$ and $ν$.

\subsubsection{Influence of $j$ on the effectivity index}\label{sec:numex:effectivity:j}
\begin{figure}
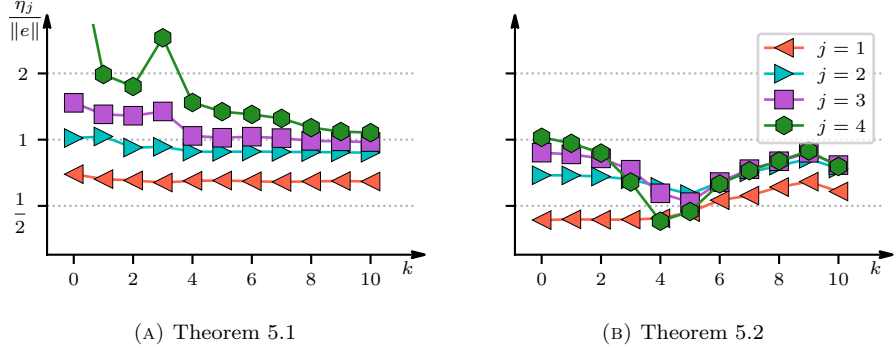

    \centering
    \subfloat[\label{fig:effectivity_j:waves2d}\cref{ex:waves}]
        {\input{figures/waves2d_estimate_j_err_V6.pgf}}
    \subfloat[\label{fig:effectivity_j:lshape}\cref{ex:lshape}]
        {\input{figures/lshape_estimate_j_err_V6.pgf}}
    \caption{
        Effectivity indices of $η_j$ with $j=1,2,3,4$, corresponding to solving the respective problems on grids $\Tri^k_\L$, obtained by \cref{alg:kl}.
        The number of V-cycles in \cref{alg:fmg} has been set to $ν=6$.
    }
    \label{fig:effectivity_j}
\end{figure}

As stated in \cref{remark:preasymptotic}, the offset $j$ must be balanced between robustness and accuracy.
While increasing $j$ theoretically improves the effectivity index, it must be chosen small enough to  yield sufficiently small $ε$ in \cref{ass:asymptotic_convergence}.
The following numerical experiments support this theory.
Estimates $η_j$ for different values of $j$ are compared against the actual error $e_\L=u-u_\L$, using the analytical solution $u$.
Since algebraic errors are not considered in \cref{thm:eta_bounds}, we make sure that the coarse grid solutions $u_\ell$ are solved up to discretization error.
This is achieved by setting $ν=6$, \ie, deliberately over solving the systems.
The resulting effectivity indices are shown in \cref{fig:effectivity_j}.

Regarding \cref{ex:waves}, the first observation in \cref{fig:effectivity_j:waves2d} is that, the larger $j$ is, the more refinement iterations $k$ are required to yield reliable estimates.
In fact, for $j>2$, the estimates are 
unusable in the first few iterations.
This can be explained by the fact that in the beginning of the adaptive procedure, the preasymptotic regime of the uniform refinement, discussed in \cref{sec:numex:waves}, is still rather large, which results in larger $ε$ in \cref{ass:asymptotic_convergence}.
As discussed in \cref{remark:preasymptotic}, this impacts $η_j$ more, as $j$ increases.

Once the coarse grid resolves the oscillations in \cref{eq:waves_solution} sufficiently well, $ε$ becomes negligibly small, hence the estimates exhibit the expected improvement for increasing $j$.
Note that, even on the final grid $\Tri^\K$, $η_4$ is strictly worse than $η_3$.
This however, is not a contradiction to the previous claim, but seems to be an artifact of the choice to stop the procedure at $\K=10$.
Apparently, for $j=4$ an even higher degree of adaptivity would be required to satisfy \cref{ass:asymptotic_convergence} with sufficiently small $ε$.

In case of the L-shape domain, even $η_4$ remains within reasonable bounds throughout the refinement procedure.
However, there is a significant drop in effectivity around $k=4$, both for $η_3$ and $η_4$, and to a lesser degree also in $η_2$.
Similar to \cref{ex:waves}, $η_1$ does not share this issue.
Instead, it appears to improve almost monotonically.
Once again, this can be explained by the asymptotic convergence of the uniform refinement, or lag thereof.
For $k=0$, uniform refinement is expected to converge asymptotically, albeit at a slower rate of $\bigO(\hmean^{\sfrac{4}{3}})$, \ie, $θ≈2^{-\sfrac{4}{3}}$ and $ε \ll 1$.
From \cref{fig:lshape_kl_l2err}, we know that convergence is also near asymptotic for $k=\K$.
Here however, the order is close to quadratic, \ie, $θ≈2^{-2}$.
Therefore, both in the beginning and the end of the adaptive procedure, $ε$ is small, hence the effectivity is improved when increasing $j$.
Somewhere in between though, there is a transition from order $\sfrac{4}{3}$ to quadratic convergence.
Consequently, in that regime, $ε$ can become as large as $(2^{2 - \sfrac{4}{3}} -1) ≈ 0.6$, which leads to a loss of effectivity, in particular for larger $j$.

While it is crucial to understand the qualitative behavior of the estimator,
verifying the numerical values of the bounds given in \cref{thm:eta_bounds}, is just as important.
Therefore, the values of $γ_j$ shown in \cref{fig:effectivity_j} are now compared to $C_1$ and $C_2$ from \cref{thm:eta_bounds}.
Here we use idealized values for the constants, assuming $ε=0$ and theoretical convergence rates of $\Oh{2}$ for \cref{ex:waves} and $\Oh{\sfrac{4}{3}}$ for \cref{ex:lshape}.
Considering the above discussion, this should be satisfied reasonably well on $\Tri^{10}$ for \cref{ex:waves} and on $\Tri^{0}$ for \cref{ex:lshape}.
Of course, in order to be useful in practice, the estimates should also be reliable in less optimal cases.
The worst case for each $j$, \ie, the effectivity index in the iteration $k$ where $|\log(\sfrac{η_j}{\|e\|})|$ is largest, is shown in \cref{table:eff_idx}.
While $η_1$ and $η_2$ are very well contained within the idealized bounds, this is not the case for $η_3$ and $η_4$.
In is important to note that these worst case values are the outliers seen in \cref{fig:effectivity_j} caused by large $ε$, and they only exceed the bounds because these are computed assuming $ε=0$.

\begin{table}
    \centering
    \newcolumntype{g}{>{\color{gray!100}}c}
    \begin{tabular}{c |c g c g c|c g c g c}
        \multicolumn{1}{c}{} &&
            \multicolumn{3}{c}{\cref{ex:waves}} &&&
            \multicolumn{3}{c}{\cref{ex:lshape}} & \\[0.5ex]
        \toprule
        $j$ && $C_1$ & $\sfrac{η_j}{\|e\|}$ & $C_2$
            &&&$C_1$ & $\sfrac{η_j}{\|e\|}$ & $C_2$  &\\ [0.5ex]
        \midrule
        1 && 0.53 & 0.64      & 1.67  &&& 0.31 & 0.43      & 2.32 & \\
        2 && 0.80 & 0.87      & 1.24  &&& 0.53 & 0.57      & 1.76 & \\
        3 && 0.93 & \bf1.47   & 1.08  &&& 0.72 & \bf0.52   & 1.37 & \\
        4 && 0.98 & \bf8.92   & 1.02  &&& 0.85 & \bf0.43   & 1.18 & \\[1ex]
        \bottomrule
    \end{tabular}
    \caption{
        Effectivity indices from \cref{fig:effectivity_j:waves2d,fig:effectivity_j:lshape} for those coarse grids $\Tri^k$ that correspond to the worst case for each offset $j$.
        Idealized constants $C_1$ and $C_2$ from \cref{thm:eta_bounds} are computed assuming theoretical convergence factors ${θ=2^{-q}}$ and $ε=0$, with $q=2$ for \cref{ex:waves} and $q=\sfrac{4}{3}$ for \cref{ex:lshape}.
        Values exceeding the bounds are highlighted in bold.
    }
    \label{table:eff_idx}
\end{table}

These results allows us to quantify the statement made in \cref{remark:preasymptotic}, regarding the balance of accuracy and robustness.
Since only $j=1,2$ are robust for all grids in our study, and $j=2$ provides the better accuracy of the both, this seems to be a good compromise.

\subsubsection{Influence of $ν$ on the effectivity index}\label{sec:numex:effectivity:v}
In this subsection, we will study the influence of algebraic errors within the iterative solution procedure
that have been excluded in the previous section. 
As previously discussed, reducing the number of V-cycles in \cref{alg:fmg} reduces increases the algebraic errors
in the coarse grid solutions $u_\ell$.
Consequently, this could have a significant impact on the effectivity of the error estimate.
On the other hand, more V-cycles means more computational work.
Therefore, $ν$, similar to $j$, must be balanced.
Here, however, the balance is between robustness and computational cost.

To investigate the sensitivity of $η_j$ with respect to the accuracy of the coarse grid solutions, the above experiment is repeated with a reduced number of V-cycles in the FMG.
For both test problems, values as small as $ν=2$ yield very similar results to those presented in \cref{sec:numex:effectivity:j}.
It is worth noting that here, even $η_2$ slightly exceeds the bounds for one of our test problems, albeit only in the third significant digit.
As above, worst case values are summarized in \cref{table:eff_idx_v2}.

\begin{table}
    \centering
    \newcolumntype{g}{>{\color{gray!100}}c}
    \begin{tabular}{c |c g c g c|c g c g c}
        \multicolumn{1}{c}{} &&
            \multicolumn{3}{c}{\cref{ex:waves}} &&&
            \multicolumn{3}{c}{\cref{ex:lshape}} & \\[0.5ex]
        \toprule
        $j$ && $C_1$ & $\sfrac{η_j}{\|e\|}$ & $C_2$
            &&&$C_1$ & $\sfrac{η_j}{\|e\|}$ & $C_2$  &\\ [0.5ex]
        \midrule
        1 && 0.53 & 0.68      & 1.67  &&& 0.31 & 0.51      & 2.32 & \\
        2 && 0.80 & \bf1.28   & 1.24  &&& 0.53 & 0.74      & 1.76 & \\
        3 && 0.93 & \bf1.67   & 1.08  &&& 0.72 & \bf0.69   & 1.37 & \\
        4 && 0.98 & \bf7.82   & 1.02  &&& 0.85 & \bf0.53   & 1.18 & \\[1ex]
        \bottomrule
    \end{tabular}
    \caption{
        Worst case effectivity indices analogous to \cref{table:eff_idx}, but with reduced solver accuracy, where only $ν=2$ V-cycles are used in the FMG.
    }
    \label{table:eff_idx_v2}
\end{table}

When reducing the number of V-cycles even further though, \ie $ν=1$, the estimates become extremely unreliable and can be considered unusable.
It is worth mentioning, that, for problems where multigrid converges at a slower rate, the effect of $ν$ on $η_j$ might be more pronounced.
However, in such cases, more V-cycles may be required irrespective of the error estimator,
in order to keep the work on the finest grid sufficiently low.


\subsubsection{Influence of $ν$ on the mesh quality}\label{sec:numex:local_indicator}
As shown above, $η_j$ becomes 
unreliable when reducing the number of V-cycles below $ν=2$.
This raises the question whether the local error indicator $η_T$ is similarly affected.
While a detailed analysis of the issue is beyond the scope of this work,
some empirical results are available , \ie, comparing the convergence of mesh sequences produced by using local estimates of different accuracy.

\begin{figure}
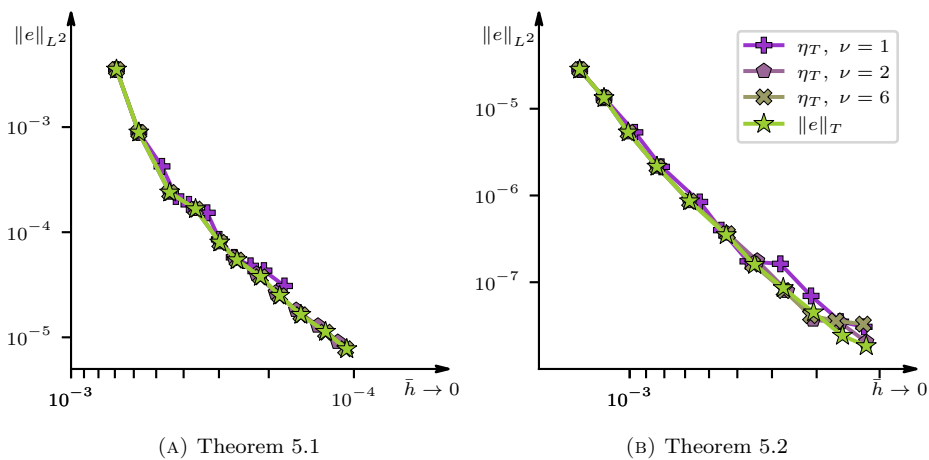

    \centering
    \subfloat[\label{fig:convergence_eta_loc:waves}\cref{ex:waves}]
    {\input{figures/waves2d_l2err_comparison_kl.pgf}}
    \subfloat[\label{fig:convergence_eta_loc:lshape}\cref{ex:lshape}]
    {\input{figures/Lshape_l2err_comparison_kl.pgf}}
    \caption{
        $L^2$ errors obtained by \cref{alg:kl}, where the AMR is based on indicators $η_T$, compared with using the actual local $L^2$ errors.
    }\label{fig:convergence_eta_loc}
\end{figure}

To that end, the experiments from \cref{sec:numex:waves,sec:numex:LShape} are repeated.
This time, however, local indicators $η_T$ are used to guide the AMR, where estimates $η_T$ based on coarse grid solutions with $ν∈\{1,2,6\}$ are considered.
The convergence of the resulting $k\ell$-refined grid sequences are shown in \cref{fig:convergence_eta_loc}.

For both problems, reducing the accuracy of the local estimator does not have a major impact on the grid convergence.
While, for \cref{ex:lshape}, the convergence rates when using $η_T$ are less smooth than with $\|e\|_T$,
the procedure appears to converge to a result of comparable quality.
However, in the final grid $\Tri^\K_\L$, the error ratio $\errratio$, as defined in \cref{eq:errratio}, increases from $1.6$ to up to $2.6$.
This is mostly due to the staircase shape of the curves in \cref{fig:convergence_eta_loc:lshape},
and improves significantly in iteration $k=9$ or $k=11$.
Overall, these results suggest that solving the coarse grid solutions to rather rough accuracy is already sufficient to yield local indicators $η_T$ that steer the $k\ell$-refinement towards a near optimal mesh sequence.

\subsection{Scalability}\label{sec:numex:scaling}
When developing and implementing algorithms for a massively parallel HPC framework such as \hyteg,
apart from convergence and accuracy, additional aspects must be considered.
For solving large practical problems, scalability is of utmost importance.
%

Therefore, a weak scaling study has been conducted, using the 3D version of the problem described in \cref{ex:waves} with $u(x,y,z) = w(x)w(y)w(z)$, and $\Omega=(0,1)^3$.
The marking strategy has been adapted such that the number of macro elements scales appropriately:
Instead of marking $10\%$ of the macro elements, $12.5\%$ are marked.
The parameters for the FMG solver are chosen based on the results discussed above, \ie,
two V-cycles on each level, with two Gauss-Seidel iterations for each pre- and post-smoothing.
On the coarse grid, CG iterations are performed until reaching a residual norm of $10^{-9}$.

In order to asses the scalability, it is sufficient to apply only a single AMR step.
That means solving on $\Tri_{\L}^0$, refining the coarse grid and solving again $\Tri_{\L}^1$.
The scheme of the study is summarized in \cref{alg:scaling}.

\begin{algorithm}[H]
    \caption{Program steps for scaling study}
    \begin{algorithmic}[1]
        \item Distribute $\Tri_0^0$ over the processes and create parallel data structures
        \item Solve for $\uh_\L$ on $\Tri_8^0$ using \textproc{FMG}
        \item Use \textproc{ErrorEstimate} to compute $η_T$
        \item Refine $\Tri_0^0$ using \cref{alg:rg3d} to obtain $\Tri_0^1$
        \item Distribute $\Tri_0^1$ over the processes and create parallel data structures
        \item Solve for $\uh_\L$ on $\Tri_8^1$ using \textproc{FMG}
    \end{algorithmic}
    \label{alg:scaling}
\end{algorithm}

The experiment is conducted on the \emph{Fritz} CPU cluster at the \emph{National High Performance Computing Center} (NHR@FAU) in Erlangen, Germany\footnote{\url{https://hpc.fau.de/systems-services/documentationinstructions/clusters/fritz-cluster}}.
Each of the 992 Nodes is equipped with 256\,GB of main memory and two \emph{Intel Xeon Platinum 8360Y}, each of which has 36 physical cores running at a clock frequency of 2.4\,GHz.
Using pure MPI parallelization, \ie, one process is started per core, leads to $n_p=72⋅n_\mathrm{nodes}$ processes.
Grid sizes and marking strategy are chosen such that $\ceil*{\sfrac{|\Tri_0^0|}{n_p}} = 3$ and $\ceil*{\sfrac{|\Tri_0^1|}{n_p}} = 7$.
The setup, together with total execution times are summarized in \cref{table:weak_scaling}.

\begin{table}
    \centering
    \begin{tabular}{lrrrrr}
        \toprule
        $n_\mathrm{nodes}$  & $n_p$  & $|\Tri_0^0|$  & $|\Tri_0^1|$  & DoF on $\Tri_8^1$  & $t_\mathrm{tot}$    \\[0.5ex]
        \midrule
         1                  &    72  &    162        &    448        &  $1.2⋅10^9$                   &  10.4\,s \\
         4                  &   288  &    750        &  1,827        &  $5.1⋅10^9$                   &  12.1\,s \\
        10                  &   720  &  2,058        &  4,695        & $13.1⋅10^9$                   &  13.4\,s \\
        21                  & 1,512  &  4,374        &  9,558        & $26.7⋅10^9$                   &  14.8\,s \\
        38                  & 2,736  &  7,986        & 17,940        & $50.1⋅10^9$                   &  16.5\,s \\
        48                  & 3,456  & 10,368        & 23,036        & $64.4⋅10^9$                   &  17.7\,s \\
        64                  & 4,608  & 13,182        & 29,114        & $81.4⋅10^9$                   &  19.7\,s \\[1ex]
        \bottomrule
    \end{tabular}
    \caption{Mesh data for the weak scaling study, as well as total execution time $t_\mathrm{tot}$ of \cref{alg:scaling}.
    }
    \label{table:weak_scaling}
\end{table}

\begin{figure}%
    \centering
    {\input{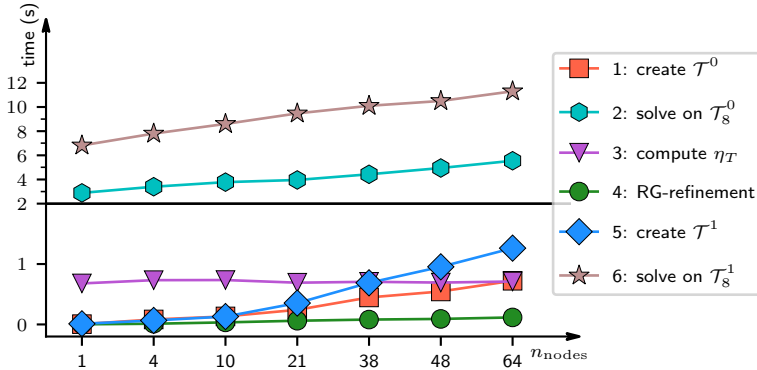}}
    \caption{
        Weak scaling results for each of the steps in \cref{alg:scaling} with the setup described in \cref{table:weak_scaling}.
        To better visualize the smaller contributions, values below 2 are shown in a different scale.
        Note that step 3 includes not only \textproc{ErrorEstimate} but also copying the data from $\uh$ to $\wh$ (line 8 in \textproc{FMG}).
    }
    \label{fig:weak_scaling}
\end{figure}

A more detailed breakdown of the compute time is given in \cref{fig:weak_scaling}.
The results are averaged over five runs of \cref{alg:scaling}.
Obviously, both RG-refinement and computing $η_T$ only take a minor fraction of the total runtime.
The latter, being a purely local operation, also exhibits near optimal weak scaling.
Given that \cref{alg:rg3d} is executed sequentially, as explained in \cref{sec:rg}, it of course does not scale in the same way.
However, the  time for executing the AMR procedure is such a minor contribution that this lack of scaling is hardly visible in the plots.
In fact, even for the largest grid in our study, the RG algorithm requires less than $1\%$ of the total time.
It is also noteworthy that, for the chosen problem -- Poisson's equation discretized with $\Pj$-elements -- solving the system is extremely cheap, so that for more complex problems the fraction of compute times spent on AMR would be even lower.
Thus, especially when considering more involved systems, such as higher order methods, curved domains, or variable coefficients,
where the solve-step becomes more expensive,
we can conclude that the sequential AMR process on the coarse grid
will not impede the overall scalability.

\begin{remark}
    We note here that the scaling of the solver is slightly suboptimal.
    This is mostly due to the coarse grid solver as employed above and the strict accuracy requirements imposed for it.
    We refer to \cite{buttari2022block, vacek2024effect} for recent studies of alternative multigrid coarse grid
    solvers and their accuracy requirements.
    Since the possible improvement is found to be of minor significance in the context of our AMR procedures,
    we refrain from a more elaborate study.
\end{remark}

  \section{Conclusion}\label{sec:conclusion}
In this work, techniques for adaptive mesh refinement on hierarchical hybrid grids were introduced.
Restricting the adaptivity to the coarse grid, the level of hierarchical structure is fully preserved,
hence the excellent scalability and performance of the associated multigrid solver is not compromised.
It was shown that this hybrid $k\ell$-refinement scheme is capable of achieving quasi optimal convergence rates, even in the presence of singularities.
Although the adaptive coarse grid refinement is carried out sequentially, we were able to demonstrate that the computational overhead is negligible compared to the cost of solving the linear system.

Furthermore, an a posteriori error estimator was presented, that is simple and extremely cost efficient, yet demonstrated to be quite effective.
The estimator is embedded into a full multigrid solver and relies on the block-wise structure of the grid.
Making maximal use of information generated automatically in a full multigrid solution process,
the computation is very cheap.
Furthermore, it can be executed fully in parallel.
In fact, the estimate itself comes at no additional cost,
but computing its norm requires a local matrix-vector product and an inner product.
The performance of the method was demonstrated in a weak scaling study, using an implementation in the HHG framework \hyteg.
Besides performance aspects, the effectivity of the error estimator was analyzed.
Analytical bounds for both local and global estimates were derived and verified in numerical studies.

Future extensions of this work will include the application to applied problems that can profit from adaptivity, such as variational inequalities
\cite{dostal2022highly} through a coupling with the PERMON library \cite{hapla2015solving}.
Furthermore we plan the implementation of an improved coarsening algorithm.
This will enable an adaptive meshing strategy that maintains a quasi-constant number of macro elements,
which is essential to better utilize the allocated compute resources in each refinement step $k$.
Another necessary extension is the interpolation between grids $\Tri^k_\ell$ and $\Tri^{k+1}_\ell$, which is currently not implemented.
Finally, we plan to study a variant of $τ$-extrapolation~\cite{tau_extrapolation, kuhn2021energy, kuhn2022implicitly, leleux2025complexity} in the context of \hyteg\ combined with the refinement strategy described in \cite{bai1987local}.
The goal here is to improve the accuracy of the approximation in \cref{eq:e_tilde} and simultaneously
exploit these techniques to obtain higher order
accuracy of the finite element approximation as well as improved error estimates without loosing the benefits of scalability.

  \section*{Acknowledgements}
The research was co-funded by the financial support of the European
Union under the REFRESH -- Research Excellence For Region Sustainability
and High-tech Industries project number CZ.10.03.01/00/22\_003/0000048
via the Operational Programme Just Transition.
This work has also received funding from the European High Performance Computing Joint Undertaking and Poland, Germany, Spain, Hungary, France, and Greece under grant agreement number: 101093457.
This publication expresses the opinions of the authors and not necessarily those of the EuroHPC JU and associated countries, which are not responsible for any use of the information contained in this publication.
The authors gratefully acknowledge the scientific support and HPC resources provided by the Erlangen National High Performance Computing Center (NHR@FAU) of the Friedrich-AlexanderUniversit\"at Erlangen-N\"urnberg (FAU).
NHR funding is provided by federal and Bavarian state authorities.
NHR@FAU hardware is partially funded by the German Research Foundation (DFG), 440719683.

  \bibliographystyle{plainnat}

\end{document}